\documentclass[11pt,a4paper]{article}

\usepackage{epic,eepic,epsfig,amsmath,amssymb,amsthm}
\usepackage{t1enc}

\usepackage{color}

\usepackage{bbm}

\def\virgp{\raise 2pt\hbox{,}}

\renewcommand{\geq}{\geqslant}
\renewcommand{\leq}{\leqslant}
\def\N{{\mathbb N}}
\def\C{\mathbb{C}}
\def\R{{\mathbb R}}

\def\virgp{\raise 2pt\hbox{,}}
\def\cdotpv{\raise 2pt\hbox{;}}

\def\1{\mathbbm{1}}
\newcommand{\ds}{\displaystyle}

\newtheorem{theorem}{Theorem}[section]

\newtheorem{proposition}[theorem]{Proposition}
\newtheorem{lemma}[theorem]{Lemma}

\newtheorem{pte}[theorem]{Property}

\theoremstyle{remark}
\newtheorem{remark}{Remark}[section]

\theoremstyle{definition}
\newtheorem{definition}{Definition}[section]

\newtheorem*{notation}{Notation}

\theoremstyle{definition}

\theoremstyle{definition}

\begin{document}

\title{\textbf{Laplacian, on the graph of the Weierstrass function}}

\author{\LARGE{\textbf{Claire David}}}

\maketitle
\centerline{Sorbonne Universit\'e}

\centerline{CNRS, UMR 7598, Laboratoire Jacques-Louis Lions, 4, place Jussieu 75005, Paris, France}

\begin{abstract}
We present, in the following, the results which enable one to build a Laplacian on the graph of the Weierstrass function, by following the approach of J. Kigami and R. S. Strichartz. Ours is made in a completely renewed framework.
\end{abstract}

\maketitle

\vskip 1cm

\vskip 1cm

\section{Introduction}

The Laplacian plays a major role in the mathematical analysis of partial differential equations. Recently, the work of J. Kigami~\cite{Kigami1989}, ~\cite{Kigami1993}, taken up by
 R.~S.~Strichartz ~\cite{Strichartz1999}, \cite{StrichartzLivre2006}, allowed the construction of an operator of the same nature, defined locally, on graphs having a fractal character: the Sierpi\'nski gasket, the Sierpi\'nski carpet, the diamond fractal, the Julia sets, the fern of Barnsley. \\

J.~Kigami starts from the definition of the Laplacian on the unit segment of the real line. For a  double-differentiable function~$u$ on~$ [0,1] $, the Laplacian~$ \Delta \, u $ is obtained as a second derivative of~$ u $ on~$ [0,1] $. For any pair~$ (u, v) $ belonging to the space of functions that are differentiable on $ [0,1] $, such that:
$$v(0)=v(1)=0$$

\noindent he puts the light on the fact that, taking into account:

$$  \displaystyle \int_0^1 \left (\Delta u \right)(x)\,v(x)\,dx=-   \displaystyle \int_0^1 u'(x)\,v'(x)\,dx
=- \displaystyle \lim_{n \to + \infty} \sum_{k=1}^n  \displaystyle \int_{ \frac{k-1}{n}}^{ \frac{k }{n}} u'(x)\,v'(x)\,dx $$

\noindent if~$\varepsilon > 0$, the continuity of~$u'$ and~$v'$ shows the existence of a natural rank~$n_0$ such that, for any integer~$n \geq n_0$, and any real number~$x$ of~$\left [\displaystyle \frac{k-1}{n}, \displaystyle\frac{k }{n} \right]$,~$1 \leq k \leq n$:

$$\left | u'  (x)- \displaystyle   \frac{ u \left (\displaystyle \frac{k}{n} \right) -u \left (\displaystyle \frac{k-1}{n} \right)}{\frac{1}{n} }\right| \leq \varepsilon
\quad , \quad \left | v'  (x)- \displaystyle   \frac{ v \left (\displaystyle \frac{k}{n} \right) -v \left (\displaystyle \frac{k-1}{n} \right)}{\frac{1}{n} } \right| \leq \varepsilon$$
\noindent the relation:

$$  \displaystyle \int_0^1 \left (\Delta u \right)(x)\,v(x)\,dx=
-\displaystyle \lim_{n \to + \infty} n \,\displaystyle \sum_{k=1}^n \left (u \left (\displaystyle \frac{k}{n} \right) -u \left (\displaystyle \frac{k-1}{n} \right) \right) \,
\left (v \left (\displaystyle \frac{k}{n} \right) -v \left (\displaystyle \frac{k-1}{n} \right) \right)
 $$

\noindent enables one to define, under a  form, the Laplacian of~$u$, while avoiding first derivatives. It thus opens the door to Laplacians on fractal domains.\\

Concretely, the weak formulation is obtained by means of Dirichlet forms, built by induction on a sequence of graphs that converges towards the considered domain. For a continuous function on this domain, its Laplacian is obtained as the renormalized limit of the sequence of graph Laplacians.  \\

If the work of J.~Kigami is, in means of analysis on fractals, seminal, it is to \textbf{Robert~S.~Strichartz} that one owes its rise. Robert~S.~Strichartz goes further than J.~Kigami : on the ground of the Sierpi\'nski gasket, he deepens, develops, exploits, generalizes, and reconstructs the classical functional spaces.\\

Strangely, the case of the \textbf{graph of the Weierstrass function}, introduced in~1872 by K.~Weierstrass~\cite{Weierstrass1872}, which presents self similarity properties, does not seem to have been considered anywhere, \textbf{neither by Robert~S.~Strichartz, neither by others}. It is yet \textbf{an obligatory passage}, in the perspective of studying {\textbf{diffusion phenomena in irregular structures.}}\\
\noindent Let us recall that, given~$\lambda \,\in\,]0,1[$, and~$b$ such that~$\lambda\,b > 1+\displaystyle \frac{3\,\pi}{2}$, the  Weierstrass function\index{Weierstrass}

$$ x \,\in\,\R \mapsto  \displaystyle \sum_{n=0}^{+\infty} \lambda^n\,\cos \left (   \pi\,b^n\,x \right)
$$

\noindent is continuous everywhere, while nowhere differentiable. The original proof, by~K.~Weierstrass~\cite{Weierstrass1872}, can also be found in~ \cite{Titschmarsh1939}. It has been completed by the one, now a classical one, in the case where~$  \lambda\, b > 1$, by~G.~Hardy \cite{Hardy1911}.


$$ $$




 After the works of A.~S.~Besicovitch and H.~D.~Ursell~\cite{BesicovitchUrsell1937}, it is Beno\^{i}t~Mandelbrot~\cite{Mandelbrot1977} who particularly highlighted the fractal properties of the graph of the Weierstrass function. He also conjectured that the Hausdorff dimension of the graph is~\mbox{$D_{\cal W}= 2 + \displaystyle \frac{\ln \lambda}{\ln b}$}. Interesting discussions in relation to this question have been given in the book of K.~Falconer \cite{Falconer1985}. A proof was given by~B.~Hunt \cite{Hunt1998} in 1998 in the case where arbitrary phases are included in each cosinusoidal term of the summation. Recently, K.~Bara\'{n}sky, B.~B\'{a}r\'{a}ny and J.~Romanowska~\cite{Baransky2014} proved that, for any value of the real number~$b$, there exists a threshold value~$\lambda_b$ belonging to the interval~\mbox{$\left ]\displaystyle \frac{1}{b},1\right [$}  such that the aforementioned dimension is equal to~\mbox{$D_{\cal W}$} for every~$b$ in~\mbox{$ \left ]\lambda_b,1 \right[$}. In \cite{Keller2017}, G.~Keller proposes what appears as a much simpler and very original proof.

\vskip 0.5cm

We have asked ourselves the following question: \textbf{given a continuous function~$u$ on the graph of the Weierstrass function, under which conditions is it possible to associate to~$u$ a function~\mbox{$ \Delta \, u $} which is, in the weak sense, its Laplacian, so that this new function~\mbox{$ \Delta \, u $} is also defined and continuous on the graph of the Weierstrass function ?}  \\

We present, in the following, the results obtained by following the approach of~J.~Kigami and R.~S.~Strichartz. Ours is made \textbf{in a completely renewed framework}, as regards, the one, affine, of the Sierpi\' nski gasket. First, we concentrate on Dirichlet forms, on the graph of the Weierstrass function, which enable us the, subject to its existence, to define the Laplacian of a continuous function on this graph. This Laplacian appears as the renormalized limit of a sequence of discrete Laplacians on a sequence of \textbf{graphs} which converge to the one of the Weierstrass function. The normalization constants related to each graph Laplacian are obtained thanks Dirichlet forms.\\

In addition to the Dirichlet forms, we have come across several delicate points: the building of a specific measure related to the graph of the function, as well as the one of spline functions on the vertices of the graph. \\

The \textbf{spectrum of the Laplacian} thus built is obtained through spectral decimation. Beautifully, as regards to the method developed by Robert S.~Strichartz in the case of the Sierpi\' nski gasket, our results come as the most natural illustration of the iterative process that gives birth to the discrete sequence of graphs.


\vskip 1cm

\section{Dirichlet forms, on the graph of the Weierstrass function}

\begin{notation}
In the following,~$\lambda$ and~$b$ are two real numbers such that:

$$0 <\lambda<1 \quad, \quad b=N_b\,\in\,\N \quad  \text{and} \quad \lambda\,N_b > 1 $$

\noindent We will consider the Weierstrass function~${\cal W}$, defined, for any real number~$x$, by:

$$   {\cal W}(x)=\displaystyle \sum_{n=0}^{+\infty} \lambda^n\,\cos \left ( 2\,  \pi\,N_b^n\,x \right)
$$

\end{notation}
\vskip 1cm

\subsection{Theoretical study}

We place ourselves, in the following, in the euclidian plane of dimension~2, referred to a direct orthonormal frame. The usual Cartesian coordinates are~$(x,y)$.\\














\vskip 1cm

\begin{pte} \textbf{Periodic properties of the Weierstrass function}\\
\noindent For any real number~$x$:

$$   {\cal W}( x+1)=\displaystyle \sum_{n=0}^{+\infty} \lambda^n\,\cos \left ( 2\,\pi\,N_b^{n }\,x +2\,\pi\,N_b^{n }\right)
=\displaystyle \sum_{n=0}^{+\infty} \lambda^n\,\cos \left ( 2\,\pi\,N_b^{n }\,x  \right)={\cal W}( x )
$$


\noindent The study of the Weierstrass function can be restricted to the interval~$[0,1[$.

\end{pte}

\vskip 1cm

\noindent By following the method developed by~J.~Kigami, we approximate the restriction~$\Gamma_{\cal W}$ to~$[0,1[ \times \R$, of the graph of the Weierstrass function, by a sequence of graphs, built through an iterative process. To this purpose, we introduce the iterated function system of the family of~$C^\infty$ maps from~$\R^2$ to~$\R^2$:
$$\left \lbrace T_{0},...,T_{N_b-1} \right \rbrace$$

\noindent where, for any integer~$i$ belonging to~$\left \lbrace 0,...,N_b-1  \right \rbrace$, and any~$(x,y)$ of~$\R^2$:
$$ T_i(x,y) =\left( \displaystyle \frac{x+i}{N_b}, \lambda\, y + \cos\left(2\,\pi \,\left(\frac{x+i}{N_b}\right)\right) \right)$$

\vskip 1cm



\vskip 1cm

\begin{lemma}

For any integer~$i$ belonging to~\mbox{$\left \lbrace  0,...,N_b-1  \right \rbrace$}, the map~$T_i$ is a bijection of~$\Gamma_{\cal W}$.

\end{lemma}

\vskip 1cm

\begin{proof}
\noindent Let~$i \,\in\,\left \lbrace  0,...,N_b-1  \right \rbrace$.
\noindent Consider a point~$\left ( y,{\cal W}(y)\right)$ of~$\Gamma_{\cal W}$, and let us look for a real number~$x$ of~$[0,1]$ such that:

$$T_i\left (x,{\cal W}(x) \right)= \left ( y,{\cal W}(y)\right)$$

\noindent One has:

$$y=\displaystyle \frac{x+i}{N_b}$$

\noindent Then:

$$x = N_b\,y-i$$

\noindent This enables one to obtain:

$${\cal W}(x)=
{\cal W}(N_b\,y-i)=\displaystyle \sum_{n=0}^{+\infty} \lambda^n\,\cos \left ( 2\,\pi\,N_b^{n+1 }\,y -2\,\pi\,N_b^{n }\,i\right)
=\displaystyle \sum_{n=0}^{+\infty} \lambda^n\,\cos \left ( 2\,\pi\,N_b^{n+1 }\,y  \right)$$

\noindent and:

$$\begin{array}{ccc}
T_i\left (x,{\cal W}(x)\right)&=&\left( \displaystyle \frac{x+i}{N_b}, \lambda\, {\cal W}(x) + \cos\left(2\,\pi \,\left(\frac{x+i}{N_b}\right)\right) \right)\\
&=&\left(y, \lambda\,\displaystyle \sum_{n=0}^{+\infty} \lambda^n\,\cos \left ( 2\,\pi\,N_b^{n+1 }\,y  \right) + \cos\left(2\,\pi \,\left(y\right)\right) \right)\\
&=&\left(y,  \displaystyle \sum_{n=0}^{+\infty} \lambda^{n+1}\,\cos \left ( 2\,\pi\,N_b^{n+1 }\,y  \right) + \cos\left(2\,\pi \,\left(y\right)\right) \right)\\
&=&\left(y,  \displaystyle \sum_{n=0}^{+\infty} \lambda^{n }\,\cos \left ( 2\,\pi\,N_b^{n  }\,y  \right)  \right)\\
&=&\left(y,{\cal W}(y)\right)
\end{array}$$

\noindent There exists thus a unique real number~$x$ in~$[0,1]$ such that:

$$T_i\left (x,{\cal W}(x) \right)= \left ( y,{\cal W}(y)\right)$$










\end{proof}

\vskip 1cm

\begin{pte}

$$   \Gamma_{\cal W} =  \underset{  i=0}{\overset{N_b-1}{\bigcup}}\,T_{i}(\Gamma_{\cal W})$$

\end{pte}

\vskip 1cm









\begin{remark}
\noindent The family~$\left \lbrace T_{0},...,T_{N_b-1} \right \rbrace$ is a family of contractions from~$\R^2$ to~$\R^2$.\\

\vskip 1cm
\begin{proof}

 \noindent For any integer~$i$ belonging to~\mbox{$\left \lbrace 0,...,N_b-1  \right \rbrace$}, we introduce the jacobian matrix of~$T_i$ such that, for any~$(x,y)$ of~$\R^2$   :

$$
DT_i(x,y)=\left (
\begin{array}{cc}
\displaystyle \frac{1}{N_b} & 0 \\
-\displaystyle \frac{2\,\pi}{N_b}\sin\left(2\,\pi\, \left( \displaystyle \frac{x+i}{N_b})\right) \right)  & \lambda
\end{array} \right )
$$

\noindent For any~$i$ of~$\left \lbrace 0,...,N_b-1  \right \rbrace$, the spectral radius of the linear map

$$\begin{array}{ccc} \R^2 & \to &\R^2 \\
(u,v)& \mapsto & DT_i(x,y) \,\left ( \begin{array}{c} u \\v \end{array} \right) \end{array}$$

\noindent  is:
 $$\displaystyle{\rho(DT_i(x,y))=\max\, \left\lbrace \displaystyle \frac{1}{N_b},\lambda\right\rbrace }= K <1$$

 \noindent Let us denote by~$\| \cdot\|_{2 }$ the euclidean norm on~$\R^2$ :

  $$\forall\, X=(x,y)\,\in\,\R^2\, : \quad \|X\|_{2 }=\|(x,y)\|_{2 }=\sqrt{x^2+y^2}$$

\noindent  For the spectral norm~$\| \cdot\|_{2,2}$, defined on the space of ~$2 \times 2$ real matrices:

 $$ A=\left (
\begin{array}{cc}
a_{11} & a_{12}\\
a_{21} & a_{22}
\end{array} \right ) \mapsto\displaystyle \sup_{ \|  X  \|_2=1}\, \left \|A\, X \right\|_2$$

 \noindent one has:

$$\forall\, (x,y) \,\in \,\R^2 \, : \quad  \| DT_i(x,y)\|_{2,2} \leq \max\, \left\lbrace \displaystyle \frac{1}{N_b},\lambda\right\rbrace  <1     $$

\noindent Thus, for any quadruplet~$(x,y,z,t)$ of real numbers:

$$
\left\| T_i(x,y)-T_i(z,t)\right\|_{2 }  \leq K \, \left\| (x,y)-(z,t) \right\|_{2 }
$$

\end{proof}

\end{remark}
\vskip 1cm

\begin{definition}
\noindent For any integer~$i$ belonging to~$\left \lbrace 0,...,N_b-1\right \rbrace  $, let us denote by:

$$P_i=(x_i,y_i)=\left(\displaystyle \frac{i}{N_b-1},\displaystyle\frac{1}{1-\lambda}\,\cos\left ( \displaystyle\frac{2\,\pi\,i}{N_b-1}\right ) \right) $$

\noindent the fixed point of the contraction~$T_i$.\\

\noindent We will denote by~$V_0$ the ordered set (according to increasing abscissa), of the points:

$$\left \lbrace P_{0},...,P_{N_b-1}\right \rbrace$$

\noindent  since, for any~$i$ of~\mbox{$\left \lbrace  0,...,N_b-2  \right \rbrace$}:

$$x_i \leq x_{i+1}$$

\noindent The set of points~$V_0$, where, for any~$i$ of~\mbox{$\left \lbrace  0,...,N_b-2  \right \rbrace$}, the point~$P_i$ is linked to the point~$P_{i+1}$, constitutes an oriented graph (according to increasing abscissa), that we will denote by~$ \Gamma_{{\cal W}_0}$.~$V_0$ is called the set of vertices of the graph~$ \Gamma_{{\cal W}_0}$.\\

\noindent For any natural integer~$m$, we set:
$$V_m =\underset{  i=0}{\overset{N_b-1}{\bigcup}}\, T_i \left (V_{m-1}\right )$$

\noindent The set of points~$V_m$, where two consecutive points are linked, is an oriented graph (according to increasing abscissa), which we will denote by~$ \Gamma_{{\cal W}_m}$.~$V_m$ is called the set of vertices of the graph~$ \Gamma_{{\cal W}_m}$. We will denote, in the following, by~${\cal N}^{\cal S}_m$ the number of vertices of the graph~$ \Gamma_{{\cal W}_m}$, and we will write:

 $$V_m = \left \lbrace {\cal S}_0^m,  {\cal S}_1^m, \hdots,  {\cal S}_{{\cal N}_m-1}^m \right \rbrace $$  

\end{definition}

\vskip 1cm

\begin{pte}
\noindent For any natural integer~$m$:
$$V_m \subset V_{m+1} $$

\end{pte}

\vskip 1cm

\begin{pte}
For any integer~$i$ belonging to~$\left \lbrace  0,...,N_b-2  \right \rbrace$:
$$ T_i\left ( P_{N_b-1} \right) =T_{i+1}\left ( P_{0} \right) $$

\end{pte}

\vskip 1cm

\begin{proof}
\noindent Since:

$$P_0=\left(0,\displaystyle\frac{1}{1-\lambda}  \right) \quad, \quad P_{N_b-1}= \left(\displaystyle \frac{ N_b-1 }{N_b-1},\displaystyle\frac{1}{1-\lambda}\,\cos\left ( \displaystyle\frac{2\,\pi\,(N_b-1)}{N_b-1}\right ) \right)
=\left( 1,\displaystyle\frac{1}{1-\lambda}\right) $$

\noindent one has :

$$ T_i\left ( P_{N_b-1} \right) =\left( \displaystyle \frac{1+i}{N_b},   \displaystyle\frac{\lambda}{1-\lambda} + \cos\left(2\,\pi \,\left(\frac{1+i}{N_b}\right)\right) \right)$$

$$ T_{i+1}\left ( P_{0} \right) =\left( \displaystyle \frac{ i+1}{N_b},   \displaystyle\frac{\lambda}{1-\lambda} + \cos\left(2\,\pi \,\left(\frac{i+1}{N_b}\right)\right) \right)$$

 \end{proof}

\vskip 1cm

\begin{pte}
\noindent The sequence~$\left ({\cal N}^{\cal S}_m \right)_{m\in\N }$ is an arithmetico-geometric one, with~${\cal N}^{\cal S}_0=N_b$ as first term:

$$\forall\,m\,\in\,\N \, : \quad {\cal N}^{\cal S}_{m+1}=N_b\, {\cal N}^{\cal S}_m- \left (N_b-2\right)$$

\noindent This leads to:

$$\forall\,m\,\in\,\N \, : \quad {\cal N}^{\cal S}_{m }=N_b^m\,\left ( {\cal N}_0-\left (N_b-2\right) \right)+ \left (N_b-2\right)
=2\, N_b^m+  N_b-2 $$

\end{pte}
\vskip 1cm

\begin{proof}

This results comes from the fact that each graph~$  \Gamma_{{\cal W}_m} $,~$m\,\in\,\N^\star$, is built from its predecessor~$  \Gamma_{{\cal W}_{m-1}} $ by applying the ~$N_b$ contractions~$T_i$,~$0 \leq i\leq N_b-1$, to the vertices of~$  \Gamma_{{\cal W}_{m-1}} $. Since, for any~$i$ of~$\left \lbrace 0,...,N_b-2  \right \rbrace$:
$$ T_i\left ( P_{N_b-1} \right) =T_{i+1}\left ( P_{0} \right) $$
\noindent the,~$N_b-2$ points appear twice if one takes into account the images of the~${\cal N}_{m-1}$ vertices of~$  \Gamma_{{\cal W}_{m-1}} $ by the whole set of contractions~$T_i$,~$0 \leq i\leq N_b-1$.

\end{proof}
\vskip 1cm

\vskip 1cm

 \begin{figure}[h!]
 \center{\psfig{height=8cm,width=9cm,angle=0,file=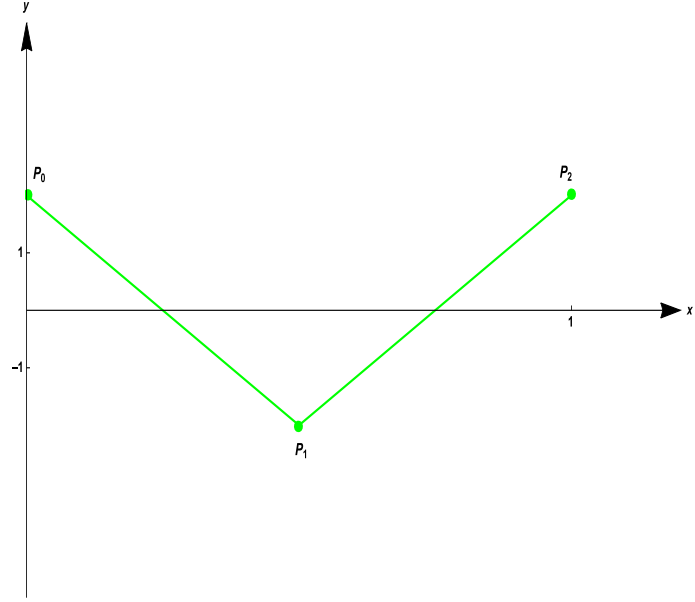}}\\
\caption{The fixed points~$P_0$,~$P_1$,~$P_2$, and the graph~$ \Gamma_{{\cal W}_0 }$, in the case where~$\lambda= \displaystyle \frac{1}{2}$, and~$N_b=3$.}
 \end{figure}

 \begin{figure}[h!]
 \center{\psfig{height=8cm,width=12cm,angle=0,file=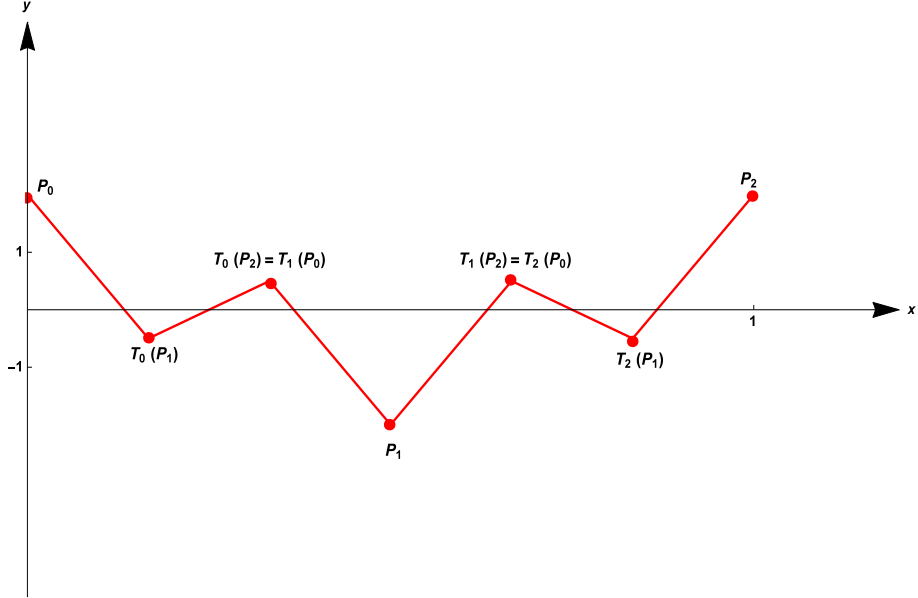}}\\
\caption{The graph~$ \Gamma_{{\cal W}_1 }$, in the case where~$\lambda= \displaystyle \frac{1}{2}$, and~$N_b=3$. $T_0(P_2)=T_1(P_0)$ et $T_1(P_2)=T_2(P_1)$.   }
 \end{figure}

 \begin{figure}[h!]
 \center{\psfig{height=10cm,width=12cm,angle=0,file=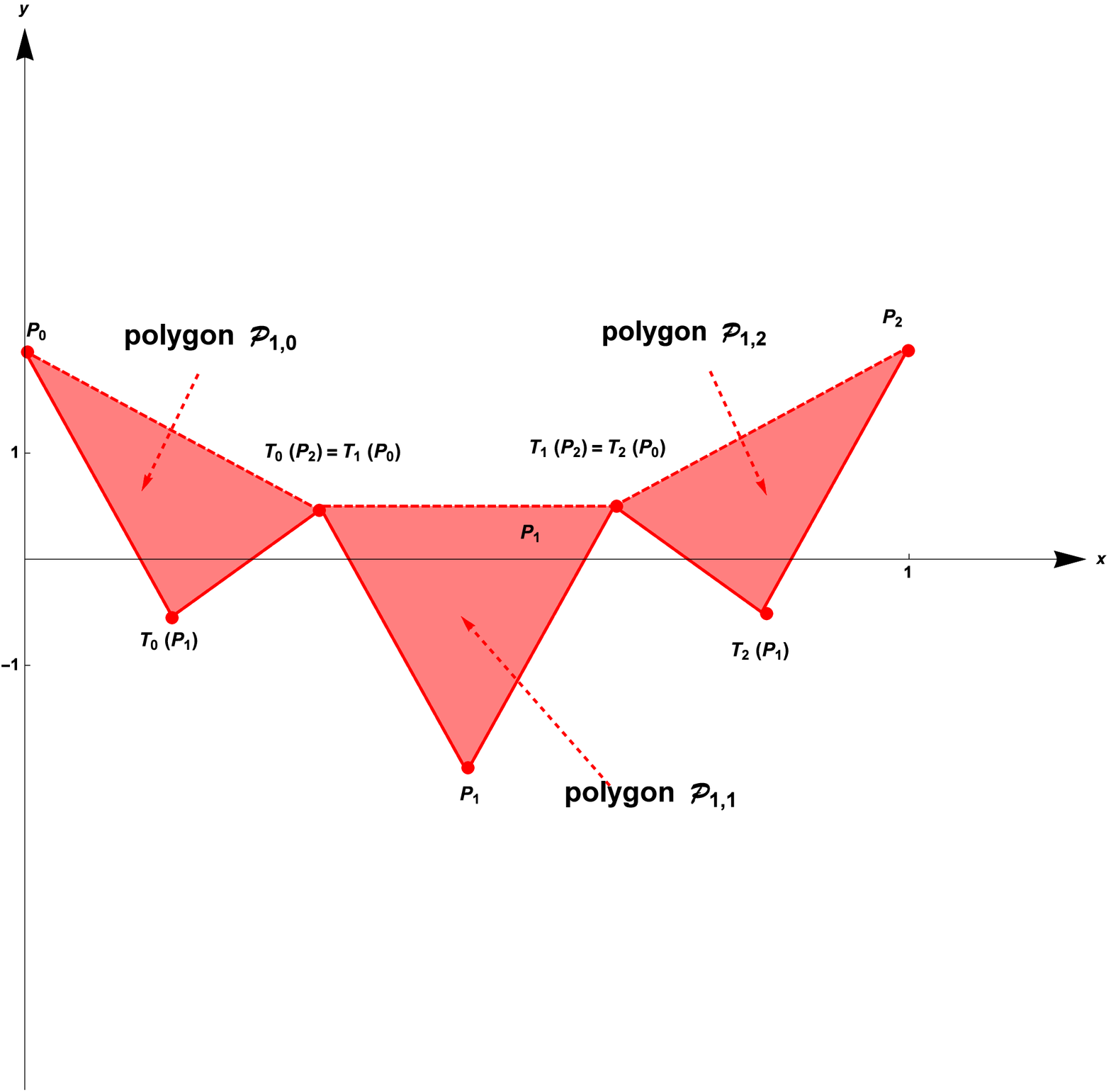}}\\
\caption{The polygons~${\cal P}_{1,0}$,~${\cal P}_{1,1}$,~${\cal P}_{1,2}$, in the case where~$\lambda= \displaystyle \frac{1}{2}$, and~$N_b=3$.}
 \end{figure}



 \begin{figure}[h!]
 \center{\psfig{height=8cm,width=12cm,angle=0,file=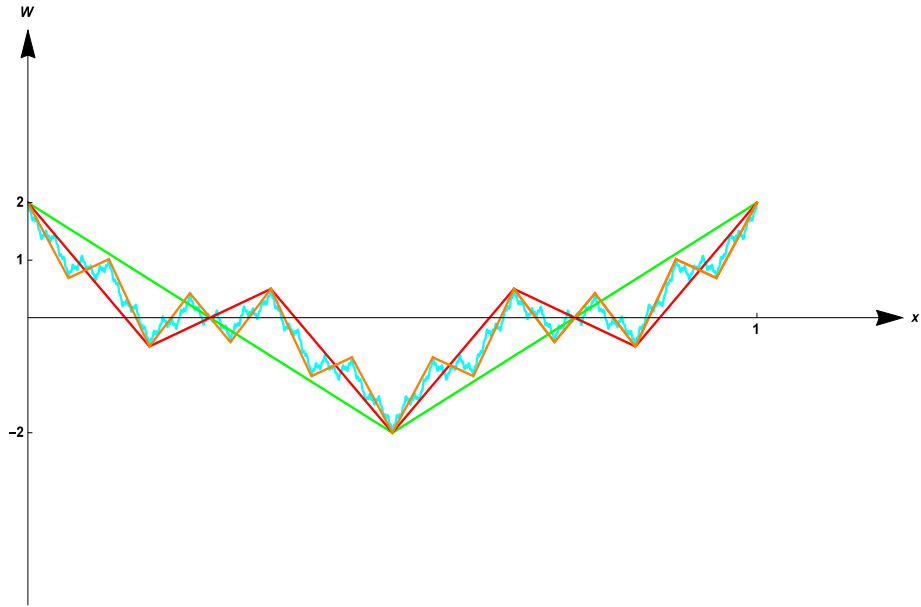}}\\
\caption{The graphs~$ \Gamma_{{\cal W}_0 }$ (in green), ~$ \Gamma_{{\cal W}_1 }$ (in red),~$ \Gamma_{{\cal W}_2 }$ (in orange),~$ \Gamma_{{\cal W} }$ (in cyan), in the case where~$\lambda= \displaystyle \frac{1}{2}$, and~$N_b=3$.    }
 \end{figure}

\vskip 1cm
\newpage

\begin{definition}\textbf{Consecutive vertices on the graph~$\Gamma_{ \cal W} $ }\\

\noindent Two points~$X$ and~$Y$ of~$\Gamma_{{ \cal W} }$ will be called \textbf{\emph{consecutive vertices}} of the graph~$\Gamma_{ \cal W} $ if there exists a natural integer~$m$, and an integer~$j $ of~\mbox{$\left \lbrace  0,...,N_b-2  \right \rbrace$}, such that:

$$X = \left (T_{i_1}\circ \hdots \circ T_{i_m}\right)(P_j) \quad \text{and} \quad Y = \left (T_{i_1}\circ \hdots \circ T_{i_m}\right)(P_{j+1})
\qquad\left \lbrace i_1,\hdots, i_m \right \rbrace \,\in\,\left \lbrace  0,...,N_b-1  \right \rbrace^m $$

\noindent or:

$$X = \left (T_{i_1}\circ  T_{i_2}\circ \hdots \circ T_{i_m}\right)\left (P_{N_b-1}\right) \quad \text{and} \quad Y =\left (T_{i_1+1}\circ T_{i_2}\hdots \circ T_{i_m} \right)(P_{0})$$

\end{definition}

\vskip 1cm

\begin{remark}

\noindent It is important to note that~$X$ and~$Y$ cannot be in the same time the images of~$P_j$ and~$P_{j+1}$,~$0 \leq j \leq N_{b-2}$, by~\mbox{$T_{i_1}\circ \hdots \circ T_{i_m}$},
~\mbox{$(i_1, \hdots,i_m)\,\in\,\left \lbrace  0,...,N_b-2  \right \rbrace$}, and of~$P_{k }$ and~$P_{k+1}$,~$0 \leq k \leq N_{b-2}$ ,\\ by~\mbox{$T_{p_1}\circ \hdots \circ T_{ p_m}$},~\mbox{$(p_1, \hdots,p_m)\,\in\,\left \lbrace  0,...,N_b-2  \right \rbrace$}. This result can be proved by induction, since, for any pair of integers~$(j,k)$ of~\mbox{$\left \lbrace  0,...,N_b-2  \right \rbrace^2$ }, for any~$i_m$ of~\mbox{$\left \lbrace  0,...,N_b-2  \right \rbrace$}, and any~$p_m$ of~\mbox{$\left \lbrace  0,...,N_b-2  \right \rbrace$}:

$$\left ( i_m \neq p_m \quad \text{and} \quad j \neq k \right)\Longrightarrow \left (T_{i_m}\left (P_j \right) \neq T_{j_m}\left (P_{k } \right) \quad \text{and} \quad T_{i_m}\left (P_j \right) \neq T_{j_m}\left (P_{k } \right) \right)$$

\noindent Each contraction~$T_i$,~\mbox{$  0 \leq i \leq N_b-1$} is indeed injective.

\noindent Since the vertices of the initial graph~$\Gamma_{{\cal W}_0}$ are distinct, one gets the expected result.

\end{remark}

\vskip 1cm

\begin{definition}
\noindent For any natural integer~$m$, the~$ {\cal N}^{\cal S}_m$ consecutive vertices of the graph~$  \Gamma_{{\cal W}_m} $ are, also, the vertices of~$N_b^m$ simple polygons~${\cal P}_{m,j}$,~\mbox{$0 \leq j \leq N_b^m-1$}, with~$N_b$ sides. For any integer~$j$ such that~\mbox{$0 \leq j \leq N_b^m-1$}, one obtains each polygon by linking the point number~$j$ to the point number~$j+1$ if~\mbox{$j = i \, \text{mod } N_b$},~\mbox{$0 \leq i \leq N_b-2$}, and the point number~$j$ to the point number~$j-N_b+1$ if~\mbox{$j =-1 \, \text{mod } N_b$}. These polygons generate a Borel set of~$\R^2$.

\end{definition}
\vskip 1cm

\begin{definition}\textbf{Polygonal domain delimited by the graph~$  \Gamma_{{\cal W}_m} $,~$m\,\in\,\N $}\\

\noindent For any natural integer~$m$, well call \textbf{polygonal domain delimited by the graph~$  \Gamma_{{\cal W}_m} $}, and denote by~\mbox{$ {\cal D} \left ( \Gamma_{{\cal W}_m}\right) $}, the reunion of the~$N_b^m$ polygons~${\cal P}_{m,j}$,~\mbox{$0 \leq j \leq N_b^m-1$}, with~$N_b$ sides.
\end{definition}
\vskip 1cm

\begin{definition}\textbf{Polygonal domain delimited by the graph~$  \Gamma_{{\cal W} } $ }\\

\noindent We will call \textbf{polygonal domain delimited by the graph~$  \Gamma_{{\cal W}} $}, and denote by~\mbox{$ {\cal D} \left ( \Gamma_{{\cal W} }\right) $}, the limit:
$$ {\cal D} \left ( \Gamma_{{\cal W} }\right)  = \displaystyle \lim_{m \to + \infty} {\cal D} \left ( \Gamma_{{\cal W}_m}\right) $$

\end{definition}

\vskip 1cm
\begin{definition}\textbf{Word, on the graph~$\Gamma_{ \cal W} $}\\

\noindent Let~$m  $ be a strictly positive integer. We will call \textbf{number-letter} any integer~${\cal M}_i$ of~\mbox{$\left \lbrace 0, \hdots, N_b-1 \right \rbrace $}, and \textbf{word of length~$|{\cal M}|=m$}, on the graph~$\Gamma_{ \cal W} $, any set of number-letters of the form:

$${\cal M}=\left ( {\cal M}_1, \hdots, {\cal M}_m\right)$$

\noindent We will write:

$$T_{\cal M}= T_{{\cal M}_1} \circ \hdots \circ  T_{{\cal M}_m}  $$

\end{definition}

\vskip 1cm

\begin{pte}

\noindent For any natural integer~$m$ :
$$\Gamma_{\cal W}=\overline{\underset{ |{\cal M}|=k \geq m  }{\overset{ }{\bigcup}}\,T_{\cal M}(\Gamma_{\cal W}) }  $$
\end{pte}
\vskip 1cm

\begin{definition}\textbf{Edge relation, on the graph~$\Gamma_{ \cal W} $}\\

\noindent Given a natural integer~$m$, two points~$X$ and~$Y$ of~$\Gamma_{{ \cal W}_m}$ will be called \emph{\textbf{adjacent}} if and only if~$X$ and~$Y$ are two consecutive vertices of~$\Gamma_{{ \cal W}_m}$. We will write:

$$X \underset{m }{\sim}  Y$$

\noindent This edge relation ensures the existence of a word~{${\cal M}=\left ( {\cal M}_1, \hdots, {\cal M}_m\right)$} of length~$ m$, such that~$X$ and~$Y$ both belong to the iterate:

$$T_{\cal M} \,V_0=\left (T_{{\cal M}_1} \circ \hdots \circ  T_{{\cal M}_m} \right) \,V_0$$

 \noindent Given two points~$X$ and~$Y$ of the graph~$\Gamma_{ \cal W} $, we will say that~$X$ and~$Y$ are \textbf{\emph{adjacent}} if and only if there exists a natural integer~$m$ such that:
$$X  \underset{m }{\sim}  Y$$
\end{definition}

\vskip 1cm

\begin{proposition}\textbf{Adresses, on the graph of the Weierstrass function}\\

\noindent Given a strictly positive integer~$m$, and a word~${\cal M}=\left ( {\cal M}_1, \hdots, {\cal M}_m\right)$ of length~$m\,\in\,\N^\star$, on the graph~$\Gamma_{ {\cal W}_m  }$, for any integer~$j$ of~$\left \lbrace 1,...,N_b-2\right \rbrace  $, any~$X=T_{\cal M}(P_j)$ of~$ V_m \setminus V_{0}$, i.e. distinct from one of the~$N_b $ fixed point~$P_i$,
  ~\mbox{$0 \leq i \leq N_b-1$}, has exactly two adjacent vertices, given by:

$$T_{\cal M}(P_{j+1})\quad \text{and} \quad T_{\cal M}(P_{j-1})$$

\noindent where:

$$T_{\cal M}  = T_{{\cal M}_1} \circ \hdots \circ  T_{{\cal M}_m}   $$

\noindent By convention, the adjacent vertices of~$T_{{\cal M} }(P_{0})  $ are~$T_{{\cal M} }(P_{1})$ and~$T_{{\cal M} }(P_{N_b-1})$,
those of~$T_{{\cal M} }(P_{N_b-1})$,~$T_{{\cal M} }(P_{N_b-2})$ and~$T_{{\cal M} }(P_{0})$ .

\end{proposition}

\vskip 1cm

\begin{pte}

The set of vertices~$\left (V_m \right)_{m \in\N}$ is dense in~$ \Gamma_{{\cal W} }$.

\end{pte}

\vskip 1cm





\begin{definition}\textbf{Power of a vertex of the graph~$ \Gamma_{{\mathcal W}_m} $,~$m\,\in\, \N^\star$}\\

	\noindent Given a strictly positive integer~$m$, a vertex~$X$ of the graph~$ \Gamma_{{\mathcal W}_m} $ will be said of power~one if~$X$ belongs to one and only one~$N_b$-gon ${\cal P}_{m,j}$,~\mbox{$0 \leq j \leq N_b^m-1$}, and of power~$\displaystyle \frac{1}{2}$ if~$X$ is a common vertex to consecutive~$N_b$-gons~${\cal P}_{m,j}$ and~${\cal P}_{m,j'}$,~\mbox{$0 \leq j \leq N_b^m-1$},~\mbox{$0 \leq j' \leq N_b^m-1$},~$j'\neq j$.\\
	In the sequel, the power of the vertex~$X$ will be denoted by:
	
	$$p(X)$$

\end{definition}

\vskip 1cm

\begin{definition}\textbf{Measure, on the domain delimited by the graph~$  \Gamma_{{\cal W} } $ }\\
	
	\noindent We will call \textbf{domain delimited by the graph~$  \Gamma_{{\cal W}} $}, and denote by~\mbox{$ {\cal D} \left ( \Gamma_{{\cal W} }\right) $}, the limit:
	$$ {\cal D} \left ( \Gamma_{{\cal W} }\right)  = \displaystyle \lim_{n \to + \infty} {\cal D} \left ( \Gamma_{{\cal W}_m}\right) $$
	
	\noindent which has to be understood in the following way: given a continuous function~$u$ on the graph~$\Gamma_{\cal W}$, and a measure with full support~$\mu$ on~$\R^2$, then:
	
	$$\displaystyle \int_{ {\cal D} \left ( \Gamma_{{\cal W} }\right)} u\,d\mu  = \displaystyle \lim_{m \to + \infty}
	\displaystyle \sum_{j=0}^{N_b^m-1}  \displaystyle \sum_{X \, \text{vertex of }{\cal P}_{m,j} }   \frac{p \left ( X \right) \,u\left ( X \right)  \,\mu \left (  {\cal P}_{m,j}  \right)}{N_b}$$
	
	\noindent We will say that~$\mu$ is a \textbf{measure, on the domain delimited by the graph~$  \Gamma_{{\cal W} } $}.
\end{definition} 
\vskip 1cm

\begin{definition}\textbf{Dirichlet form} (we refer to the paper \cite{Beurling1985}, or the book \cite{Fukushima1994})\\

 \noindent Given a measured space~$(E, \mu)$, a \emph{\textbf{Dirichlet form}} on~$E$ is a bilinear symmetric form, that we will denote by~${\cal E}$,
 defined on a vectorial subspace~$D$ dense in $L^2(E, \mu) $, such that:\\

\begin{enumerate}

\item For any real-valued function~$u$ defined on~$D$ :  ${\cal E}(u,u) \geq 0$.

\item   $D$, equipped with the inner product which, to any pair~$(u,v)$ of~$D \times D $, associates:

 $$  (u,v)_{\cal{E}}  = (u,v)_{L^2(E,\mu)} + {\cal {E}}(u,v)$$

is a Hilbert space.

\item For any real-valued function~$u$ defined on~$D$, if:
$$ u_\star = \min\, (\max(u, 0) , 1) \,\in \,D$$

\noindent then : ${ \cal{E}}(u_\star,u_\star)\leq { \cal{E}}(u,u)$ (Markov property, or lack of memory property).

\end{enumerate}

\end{definition}

\vskip 1cm

\begin{definition}\textbf{Dirichlet form, on a finite set} (\cite{Kigami2003})\\

 \noindent Let~$V$ denote a finite set~$V$, equipped with the usual inner product which, to any pair~$(u,v)$ of functions defined on~$V$, associates:

  $$(u,v)= \displaystyle \sum_{p\in  V} u(p)\,v(p)$$

  \noindent A \emph{\textbf{Dirichlet form}} on~$V$ is a symmetric bilinear form~${\cal E}$, such that:\\

\begin{enumerate}

\item For any real valued function~$u$ defined on~$V$:  ${\cal E}(u,u) \geq 0$.

\item   $  {\cal {E}}(u,u)= 0$ if and only if~$u$ is constant on~$V$.

\item For any real-valued function~$u$ defined on~$V$, if:
$$ u_\star = \min\, (\max(u, 0) , 1)  $$

\noindent i.e. :

$$\forall \,p \,\in\,V \, : \quad u_\star(p)= \left \lbrace \begin{array}{ccc} 1 & \text{if}& u(p) \geq 1 \\u(p) & \text{si}& 0 <u(p) < 1 \\0  & \text{if}& u(p) \leq 0 \end{array} \right.$$

\noindent then: ${ \cal{E}}(u_\star,u_\star)\leq { \cal{E}}(u,u)$ (Markov property).

\end{enumerate}

\end{definition}

\vskip 1cm

\begin{remark}

In order to understand the underlying theory of Dirichlet forms, one can only refer to the work of~A.~Beurling and J.~Deny~\cite{Beurling1985}. The Dirichlet space~$\cal D$ of fonctions~$u$, complex valued functions, infinitely differentiable, the support of which belongs to a domain~\mbox{$\omega \subset \R^p$},~\mbox{$p \,\in\,\N^\star$}, is equipped with the hilbertian norm:

$$ u \mapsto \|u\|_{\cal D}= \displaystyle \int_\omega |\text{grad}\, u(x) |^2\,dx$$

\noindent If the complement set of~$\omega$ is not "too small", the space~$\cal D$ can be completed by adding functions defined almost everywhere in~$\omega$. The space thus obtained~${\cal D}_\omega$, equipped with the Lebesgue measure~$\xi$, satisfies the following properties:

\begin{enumerate}
\item[\emph{i}.] For any compact~$K \subset \omega$, there exists a positive constant~$C_K$ such that, for any~$u$ of~${\cal D}_\omega$:

$$\displaystyle \int_K |u(x)|\,d\xi(x) \leq C_K \, \|u\|_{{\cal D}_\omega}$$

\item[\emph{ii}.] If one denotes by~$\cal C$ the space of complex-valued, continuous functions with compact support, then:~\mbox{${\cal C } \cap {\cal D}_\omega$} is dense in~$\cal C$ and in~${\cal D}_\omega$.

\item[\emph{iii}.] For any contraction of the complex plane, and any~$u$ of~${\cal D}_\omega$:

$$T\,u \,\in\, {\cal D}_\omega \quad \text{and} \quad \|T\,u\|_{{\cal D}_\omega} \leq  \|u\|_{{\cal D}_\omega} $$

\end{enumerate}

\noindent The Dirichlet space~${\cal D}_\omega$ is generated by the Green potentials of finite energy, which are defined in a direct way, as the functions~$u$ of~${\cal D}_\omega$ such that there exists a Radon measure~$\mu$ such that:

$$\forall \, \varphi \, \in\, {\cal C } \cap {\cal D}_\omega\, : \quad \left (u, \varphi \right)=\displaystyle \int_\omega   \bar{\varphi}\,d\mu$$
\noindent Such a map~$u$ will be called \textbf{potential generated by~$\mu$}.

\noindent The linear map~$\Delta$ which, to any potential~$u$ of~${\cal D}_\omega$, associates the measure~$\mu$ that generates this potential, is called \textbf{generalized Laplacian for the space~${\cal D}$}.\\

\noindent It is interesting to note that the original theory of Dirichlet spaces concerned functions defined on a Hausdorff space (separated espace 
), with a positive Radon measure of full support (every non-empty open set has a strictly positive measure).

\end{remark}

\vskip 1cm

\begin{remark} One may wonder why the Markov property is of such importance in our building of a Laplacian ? Very simply, the lack of memory - or the fact that the future state which corresponds, for any natural integer~$m$, to the values of the considered function on the graph~$\Gamma_{{\cal W}_{m+1}}$, depends only of the present state, i.e. the values of the function on the graph~$\Gamma_{{\cal W}_{m }}$, accounts for the building of the Laplacian step by step.

\end{remark}

\vskip 1cm

\vskip 1cm

\begin{theorem}\textbf{An upper bound and a lower bound, for the box-dimension of the graph~$\Gamma_{\mathcal W}$}\\
	
	For any integer~$j$ belonging to~\mbox{$\left \lbrace 0, 1, \hdots, N_b-2 \right \rbrace $}, each natural integer~$m$, and each word~$\mathcal M$ of length~$m$, let us consider the rectangle, whose sides are parallel to the horizontal and vertical axes, of width:

	$$L_m= x\left ( T_{\mathcal M}  \left (P_{j+1 } \right) \right)-x\left ( T_{\mathcal M}  \left (P_{j } \right)\right)=
	\displaystyle \frac{1}{(N_b-1)\, N_b^m} $$

	\noindent	 and height~$|h_{j,m}|$, such that the points~\mbox{$T_{\mathcal M}  \left (P_{j  }\right) $} and~\mbox{$T_{\mathcal M}  \left (P_{j+1 } \right )$} are two vertices of this rectangle.\\
	\noindent	 Then:\\
	
	\noindent i. When the integer~$N_b$ is odd:
	
	\footnotesize
	$$L_m^{2-D_{\mathcal W}}  \,(N_b-1)^{2-D_{\mathcal W}}\,  \left \lbrace   
	\displaystyle\frac{2}{1-\lambda}\,    \sin\left ( \displaystyle\frac{ \pi }{ N_b-1 }\right )  \displaystyle \min_{0 \leq j \leq N_b-1 }  \left | \sin\left ( \displaystyle\frac{ \pi\,(2\,j+1 )}{  N_b-1 }\right )   \right| - \displaystyle  \frac{ 2\,\pi   }{  N_b \,(N_b-1)}  \,
	\displaystyle  \frac{ 1 }{   \lambda \, N_b  -1}  \right \rbrace  \leq |h_{j,m}|  $$
	
	\normalsize

	\noindent ii. When the integer~$N_b$ is even:

	\tiny
	$$L_m^{2-D_{\mathcal W}}   (N_b-1)^{2-D_{\mathcal W}}   \max  \left \lbrace   \! \! 
	\displaystyle\frac{2}{1-\lambda}     \sin\left ( \displaystyle\frac{ \pi }{ N_b-1 }\right ) \! \!  
	\displaystyle \min_{0 \leq j \leq N_b-1 } \! \!  \left | \sin\left ( \displaystyle\frac{ \pi\,(2\,j+1 )}{  N_b-1 }\right )   \right| \! \! -\! \!  \displaystyle  \frac{ 2\,\pi   }{  N_b \,(N_b-1)}   
	\displaystyle  \frac{ 1 }{   \lambda \, N_b  -1} , 
	\displaystyle  \frac{4   }{  N_b^{2}  }      \displaystyle  \frac{1-  N_b^{ -2}  }
	{    N_b^2-1  }  
	\right \rbrace  \! \! \leq \!  |h_{j,m}| \!  $$
	
	\normalsize
	
	\noindent	 Also:
	
	$$ |h_{j,m}| \leq  \eta_{ {\mathcal W}   }\,L_m^{2-D_{\mathcal W}} \, (N_b-1)^{2-D_{\mathcal W}}\,  $$
	
	\noindent where the real constant~\mbox{$  \eta_{ {\mathcal W}   }$} is given by :
	
	$$ \eta_{ {\mathcal W}   }  = 2\, \pi^2\,  \left \lbrace
	\displaystyle  \frac{   (2\,N_b-1)\, \lambda\, (N_b^2-1) } {(N_b-1)^2 \, (1- \lambda )\,(\lambda \,N_b^{   2}-1) }    +
	\displaystyle  \frac{  2\, N_b } {   (\lambda \,N_b^{ 2  }-1)\, (\lambda \,N_b^{ 3  }-1)  } \right \rbrace .
	$$

\end{theorem}

\vskip 1cm

\begin{notation}
	
	In the sequel, we set, for any natural integer~$m$:

	$$ h_m = L_m^{2-D_{\mathcal W}}=\displaystyle   \frac{N_b^{( D_{\mathcal W}-2 )\,m}}{(N_b-1)^{2-D_{\mathcal W}}} 
	$$
	
	\noindent and:

	$$ h  =  N_b^{( D_{\mathcal W}-2 ) }
	$$
	

	
	\end{notation}



	
		\vskip 1cm

\begin{definition}\textbf{Energy, on the graph~$\Gamma_{ {\cal W }_m }$,~$m \,\in\,\N$, of a pair of functions}\\

 \noindent Let~$m$ be a natural integer, and~$u$ and~$v$ two real valued functions, defined on the set

 $$V_m = \left \lbrace {\cal S}_0^m,  {\cal S}_1^m, \hdots,  {\cal S}_{{\cal N}_m^{\cal S}-1}^m \right \rbrace $$

 \noindent of the~${\cal N}_m^{\cal S}$ vertices of~$\Gamma_{ {\cal W}_m  }$.\\






 \noindent It appears as natural to introduce \textbf{the energy, on the graph~$\Gamma_{ {\cal W }_m }$, of the pair of functions~$(u,v)$}, as:

$$\begin{array}{ccc}
  {\cal{E}}_{\Gamma_{{ \cal W }_m} }(u,v)
 &= & \displaystyle \sum_{i=0}^{{\cal N}_m^{\cal S}-2}  \left (\displaystyle \frac{u \left ({\cal S}_i^m \right)-u \left ({\cal S}_{i+1}^m \right)}{h_m}\right )\,
 \left (\displaystyle \frac{v \left ({\cal S}_{i }^m \right)-v \left ({\cal S}_{i+1}^m \right)}{h_m} \right )
 \end{array}
$$

 \noindent For the sake of simplicity, we will write it under the form:

$$ {\cal{E}}_{\Gamma_{{ \cal W }_m} }(u,v)=\displaystyle \frac{1}{h_m^2}\,\displaystyle \sum_{X  \underset{m }{\sim}  Y}\left (u(X)-u(Y)\right )\,\left(v(X)-v(Y)\right) $$

\end{definition}

\vskip 1cm

\begin{pte}

 \noindent Given a natural integer~$m$, and a real-valued function~$u$, defined on the set of vertices of~$\Gamma_{ {\cal W}_m} $, the map, which, to any pair of real-valued, continuous functions~$(u,v)$ defined on the set~$V_m $ of the~${\cal N}_m$ vertices of~$\Gamma_{ {\cal W}_m  }$, associates:
$$ {\cal{E}}_{\Gamma_{{\cal W}_m}}(u,v)=\displaystyle \frac{1}{h_m^2}\,\displaystyle \sum_{X  \underset{m }{\sim}  Y}\left (u(X)-u(Y)\right )\, \left (v(X)-v(Y)\right ) $$

\noindent is a Dirichlet form on~$\Gamma_{ {\cal W }_m }$.\\
\noindent Moreover:

$${\cal{E}}_{\Gamma_{{\cal W}_m}  }(u,u)=0 \Leftrightarrow u  \text{ is constant}$$

\end{pte}

\vskip 1cm

\begin{proposition}\textbf{Harmonic extension of a function, on the graph of the Weierstrass function}\\

\noindent For any strictly positive integer~$m$, if~$u$ is a real-valued function defined on~$V_{m-1}$, its \textbf{harmonic extension}, denoted by~$ \tilde{u}$, is obtained as the extension of~$u$ to~$V_m$ which minimizes the energy:

$$  {\cal{E}}_{\Gamma_{{\cal W}_m}}(\tilde{u},\tilde{u})=\displaystyle \sum_{X \underset{m }{\sim} Y} \displaystyle \frac{(\tilde{u}(X)-\tilde{u}(Y))^2}{h_m^2} $$

\noindent The link between~$   {\cal{E}}_{\Gamma_{{\cal W}_m}}$ and~$  {\cal{E}}_{\Gamma_{{\cal W}_{m-1}}}$ is obtained through the introduction of two strictly positive constants~$r_m$ and~$r_{m-1}$ such that:

$$  r_{m }\, \displaystyle \sum_{X \underset{m  }{\sim} Y} \displaystyle \frac{(\tilde{u}(X)-\tilde{u}(Y))^2}{h_m} =  \displaystyle \frac{r_{m-1}}{h^2}\,\displaystyle  \sum_{X \underset{m-1 }{\sim} Y} \displaystyle \frac{(u(X)-u(Y))^2}{h_{m-1}^2}$$

\noindent In particular:

$$  r_{1 }\, \displaystyle \sum_{X \underset{1  }{\sim} Y} (\tilde{u}(X)-\tilde{u}(Y))^2 =  r_{0}\,\displaystyle  \sum_{X \underset{0 }{\sim} Y} (u(X)-u(Y))^2$$

\noindent For the sake of simplicity, we will fix the value of the initial constant:~$r_0=1$. One has then:

$$ {\cal{E}}_{\Gamma_{{\cal W}_1}}(\tilde{u},\tilde{u})= \displaystyle \frac{1}{ r_{1 }\,h^2}\,  {\cal{E}}_{\Gamma_{{\cal W}_0}}(\tilde{u},\tilde{u})$$

\noindent Let us set:

$$r = \displaystyle \frac{1}{r_{1 }} $$



\noindent and:

$$  {\cal{E}}_{m}(u)=r_m\, \sum_{X \underset{m }{\sim} Y} \displaystyle \frac{(\tilde{u}(X)-\tilde{u}(Y))^2}{h_m^2} $$

\noindent Since the determination of the harmonic extension of a function appears to be a local problem, on the graph~$\Gamma_{{\cal W}_{ m-1}}$, which is linked to the graph~$\Gamma_{{\cal W}_{ m }}$ by a similar process as the one that links~$\Gamma_{{\cal W}_{  1}}$ to~$\Gamma_{{\cal W}_{ 0}}$, one deduces, for any strictly positive integer~$m$:

$$ {\cal{E}}_{\Gamma_{{\cal W}_m}}(\tilde{u},\tilde{u})= \displaystyle \frac{1}{ r_{1 }}\,  {\cal{E}}_{\Gamma_{{\cal W}_{m-1}}}(\tilde{u},\tilde{u})$$

\noindent By induction, one gets:

$$r_m=r_1^m\,r_0=r^{-m} =N_b^m$$


\noindent If~$v$ is a real-valued function, defined on~$V_{m-1}$, of harmonic extension~$ \tilde{v}$, we will write:

$$  {\cal{E}}_{m}(u,v)=r^{-m}\, \sum_{X \underset{m }{\sim} Y} 
\displaystyle \frac{(\tilde{u}(X)-\tilde{u}(Y)) \, (\tilde{v}(X)-\tilde{v}(Y))}{h_m^2} $$

\noindent For further precision on the construction and existence of harmonic extensions, we refer to~\emph{\cite{Sabot1987}}.
\end{proposition}

\vskip 1cm

\noindent  \textbf{Nota Bene:}\\

  The above latter energy formula writes:

$$  {\cal{E}}_{m}(u,u)=N_b^{(5-2\, D_{\mathcal W} )\, m}\, \sum_{X \underset{m }{\sim} Y} 
 \left (\tilde{u}(X)-\tilde{u}(Y)\right )^2 $$

  We would like to lay the emphazis upon the fact that this work involves a mixt approach, using both the methods of~J.~Kigami and~R.~S.~Strichartz, and the one by~U.~Mosco for fractal curves~\cite{UmbertoMosco2002}, which takes into account topology and geometry, by the means of a quasi-distance, built from the eucidean one~$d_{eucl}$ between adjacent points~$X$ and~$Y$ such that~${X \underset{m }{\sim} Y} $ :

$$d(X,Y)= \left ( d_{eucl} (X,Y)\right)^\delta$$

\noindent where~$\delta$ is a real constant. The related energy writes:

$$  {\cal{E}}_{m}(u)=  \sum_{X \underset{m }{\sim} Y} \displaystyle \frac{(\tilde{u}(X)-\tilde{u}(Y))^2}{d^2(X,Y)} $$

   Yet, one cannot apply this method in our case, the constant~$ \delta $ being a priori determined, either by decimation, either in order  to verify the Gaussian principle. It is impossible to determine this constant in a non-affine framework like that of the graph of the Weierstrass function.  \\

The question one may ask is wether one may be sure that~$N_b^{(5-2\, D_{\mathcal W} )\, m}$ is the right constant in the present case ? The question is not an inocuous one, in so far as the value of this constant directly affects the spectra of the related Laplacian.\\

The point is that, by construction, our energies satisfy the maximum principle. Also, the value of this constant joins the one at stake in the value of the population spectral density of fractional Brown functions evoked in~\cite{Mandelbrot1982}. It is obtained The links between randomized forms of the Weierstrass functions and fractional Brownian motion incline us to think that we are in the good direction.


\vskip 1cm







\begin{notation}
	\label{Am}

\noindent Given a strictly positive integer~$m$, let us consider a vertex~$X$ of the graph~${\Gamma}_{{\cal W}_m}$.  Two configurations can occur:\\

\begin{enumerate}
	
	\item[\emph{i}.] the vertex~$X$ belongs to one and only one polygon with~$N_b$ sides,~${\cal P}_{m,j}$,~\mbox{$0 \leq j \leq N_b^m-1$}. We set:
	
	$${\cal A}_m= \mu \left ( {\cal P}_{m,j}  \right )$$

	\item[\emph{ii}.] the vertex~$X$ is the intersection point of two polygons with~$N_b$ sides,~${\cal P}_{m,j}$ and~${\cal P}_{m,j+1}$,~\mbox{$0 \leq j \leq N_b^m-2$}. We set:

	$${\cal A}_m=\displaystyle \frac {1}{2}\,\left \lbrace  \mu \left ( {\cal P}_{m,j}  \right )+ \mu \left ( {\cal P}_{m,j+1} \right )\right  \rbrace $$
\end{enumerate}

\end{notation}

\vskip 1cm

\begin{definition}\textbf{Dirichlet form, for a pair of continuous functions defined on the graph~$\Gamma_{ {\cal W }  }$}\\

 \noindent We define the Dirichlet form~$\cal{E}$ which, to any pair of real-valued, continuous functions~$(u,v)$ defined on the graph~$\Gamma_{ {\cal W }  }$, associates, subject to its existence:

$$
  {\cal{E}} (u,v)=   
  \displaystyle \lim_{m \to + \infty}\displaystyle \sum_{X  \underset{m }{\sim}  Y} r^{-m}\, 
  \displaystyle\frac{ \left (u_{\mid V_m}(X)-u_{\mid V_m}(Y)\right )\,\left(v_{\mid V_m}(X)-v_{\mid V_m}(Y)\right)}
  {h_m^2}\, \frac{{\cal A}_m}{N_b}$$

\end{definition}

\vskip 1cm

\begin{definition}\textbf{Normalized energy, for a continuous function~$u$, defined on the graph~$\Gamma_{ {\cal W }  }$}\\
\noindent Taking into account that the sequence~$\left (\mathcal{E}_m\left ( u_{\mid V_m} \right)\right)_{m\in\N}$ is defined on
$$\ds{V_\star =\underset{{i\in \N}}\bigcup \,V_i}$$

\noindent one defines the \textbf{normalized energy}, for a continuous function~$u$, defined on the graph~$\Gamma_{ {\cal W }  }$, by:

$$
{\cal{E}} (u)=   
\displaystyle \lim_{m \to + \infty}\displaystyle \sum_{X  \underset{m }{\sim}  Y} r^{-m}\, 
\displaystyle\frac{ \left (u_{\mid V_m}(X)-u_{\mid V_m}(Y)\right )^2}
{h_m^2}\, \frac{{\cal A}_m}{N_b}$$

\end{definition}

\vskip 1cm

\vskip 1cm

\begin{notation}
\noindent We will denote by~$\text{dom}\,{\cal E}$ the subspace of continuous functions defined on~$\Gamma_{\cal W}$, such that:

$$\mathcal{E}(u)< + \infty$$

\end{notation}

\vskip 1cm

\begin{notation}
\noindent We will denote by~$\text{dom}_0\,{ \cal E}$ the subspace of continuous functions defined on~$\Gamma_{\cal W}$, which take the value on~$V_0$, such that:

$$\mathcal{E}(u)< + \infty$$

\end{notation}

\vskip 1cm




\vskip 1cm




\section{Laplacian of a continuous function, on the graph of the Weierstrass function}
\label{sec:Laplacian on the graph of the Weierstrass function}

\subsection{Theoretical aspect}







\vskip 1cm
\begin{pte}\textbf{Building of a specific measure, for the domain delimited by the graph of the Weierstrass function}\\
	
	\noindent The Dirichlet forms mentioned in the above require a positive Radon measure with full support. In auto-affine configurations (the Sierpi\'nski Gasket for instance), the choice of a self-similar measure, which is, mots of the time, built with regards to a reference set, of measure~1, appears, first, as very natural. More generally, R.~S.~Strichartz~\cite{Strichartz1995},~\cite{Strichartz1999}, showed that one can simply consider auto-replicant measures~$ \tilde{\mu}$, i.e. measures~$ \mu$ such that:

	$$ \mu= \displaystyle \sum_{i=0}^{N_b-1} \mu_i\,\mu\circ T_i^{-1} \qquad (\star)$$
	
	\noindent where~$\left (\mu_i\right)_{0 \leq i \leq N_b-1}$ denotes a family of strictly positive pounds.\\
	

	\noindent The non-affine framework makes it clear that there cannot exist such constant coefficients~$\mu_i$. It appears more realistic that they depend on the order~$m\,\in\,\N^\star$ of the iteration, as:

	$$ \mu = \displaystyle \sum_{i=0}^{N_b-1} \mu_{m,i}\,\mu \circ T_i^{-1} \qquad (\star_m)$$

	\noindent Let us thus denote by~$  \mu_{\cal L}$ the Lebesgue measure on~$\R^2$, and start with the normalized measure:
	
	$$  \mu  = \displaystyle \frac{\mu_{\cal L}}{\mu_{\cal L} \left ( {\cal P}_0\right)}$$

	\noindent and look for a family of strictly positive pounds~$\left (\mu_{1,i}\right)_{0 \leq i \leq N_b-1}$ such that:

	$$\mu   \left ( T_0\left ({\cal P}_0\right)  \cup  T_1\left ({\cal P}_0\right)  \cup  T_2\left ({\cal P}_0\right)\cup \hdots  \cup  T_{N_b-1}\left ({\cal P}_0\right) \right)
	=\displaystyle \sum_{i=0}^{N_b-1} \mu_{1,i}\,  {\mu} \left (  {\cal P}_0\right)=\displaystyle \sum_{i=0}^{N_b-1} \mu_{1,i}$$

	\noindent The convenient choice, for any~$i$ of~\mbox{$ \left \lbrace 0, \hdots, N_b-1 \right \rbrace$}, is:
	$$\mu_{1,i}=  \mu \left ( T_i\left ({\cal P}_0\right) \right)$$

	\noindent Now, given a strictly positive integer~$m$,  let us look for a family of strictly positive pounds~$\left (\mu_{1,i}\right)_{0 \leq i \leq N_b-1}$ such that:

	$$ \mu = \displaystyle \sum_{i=0}^{N_b-1} \mu_{m,i}\,\mu_0\circ T_i^{-1} \qquad (\star_m)$$
	
	\noindent where~$\left (\mu_{m,i}\right)_{0 \leq i \leq N_b-1}$ denotes a family of strictly positive pounds.\\
	
	\noindent Relation~$(\star_m)$ yields, for any set of polygons~${\cal P}_{m,j}$,~\mbox{$m\,\in\,\N$},~\mbox{$0 \leq j \leq N_b^m-1$}, with~$N_b$ sides:

	$$ \mu  \left ( \underset{0 \leq j \leq N_b^m-1 }{\bigcup} {\cal P}_{m,j} \right) =   \displaystyle \sum_{i=0}^{N_b-1} \mu_{m,i}\, \mu \left (    \left ( \underset{0 \leq j \leq N_b^{m-1}-1 }{\bigcup}{\cal P}_{m-1,j} \right) \right)$$
	
	\noindent and, in particular:
	
	$$\mu  \left ( T_0\left ( \underset{0 \leq j \leq N_b^{m-1}-1 }{\bigcup}{\cal P}_{m-1,j} \right)  \cup      \hdots  \cup  T_{N_b-1}\left ( \underset{0 \leq j \leq N_b^{m-1}-1 }	{\bigcup}{\cal P}_{m-1,j} \right) \right)
	=\displaystyle \sum_{i=0}^{N_b-1} \mu_{m,i}\,\mu\left (   \underset{0 \leq j \leq N_b^{m-1}-1 }	{\bigcup}{\cal P}_{m-1,j} \right)$$
	
	\noindent i.e.:

	$$\displaystyle \sum_{i=0}^{N_b-1}  \mu  \left (  T_i\left ( \underset{0 \leq j \leq N_b^{m-1}-1 }{\bigcup}{\cal P}_{m-1,j}\right)\right)
	=\displaystyle \sum_{i=0}^{N_b-1} \mu_i \,\mu \left (   \underset{0 \leq j \leq N_b^{m-1}-1 }{\bigcup}{\cal P}_{m-1,j}\right)$$

	\noindent The convenient choice, for any~$i$ of~\mbox{$ \left \lbrace 0, \hdots, N_b-1 \right \rbrace$ }, is:
	$$\mu_{m,i}= \displaystyle \frac{  \mu_0\left ( T_i\left ( \underset{0 \leq j \leq N_b^{m-1}-1 }{\bigcup}{\cal P}_{m-1,j}\right) \right) }{\mu  \left (   \underset{0 \leq j \leq N_b^{m-1}-1 }{\bigcup}{\cal P}_{m-1,j}\right)}$$

\end{pte}

\vskip 1cm

\begin{remark}
	
	The above result appears as as interesting splitting, fitted for the polygonal domain~\mbox{$ {\cal D} \left ( \Gamma_{{\cal W} }\right) $}.  
\end{remark}

\vskip 1cm

\begin{definition}\textbf{Laplacian of order~$m\,\in\,\N^\star$}\\

\noindent For any strictly positive integer~$m$, and any real-valued function~$u$, defined on the set~$V_m$ of the vertices of the graph~$\Gamma_{{\cal W}_m}$, we introduce the Laplacian of order~$m$,~$\Delta_m(u)$, by:

$$\Delta_m u(X) = \displaystyle \frac{1}{h_m^2}\, \displaystyle\sum_{Y \in V_m,\,Y\underset{m}{\sim} X}
 \left (u(Y)-u(X)\right)  \quad \forall\, X\,\in\, V_m\setminus V_0 $$

\end{definition}

\vskip 1cm

\begin{definition}\textbf{Harmonic function of order~$m\,\in\,\N^\star$}\\

\noindent Let~$m$ be a strictly positive integer. A real-valued function~$u$,defined on the set~$V_m$ of the vertices of the graph~$\Gamma_{{\cal W}_m}$, will be said to be \textbf{harmonic of order~$m$} if its Laplacian of order~$m$ is null:

$$\Delta_m u(X) =0 \quad \forall\, X\,\in\, V_m\setminus V_0 $$

\end{definition}

\vskip 1cm

\begin{definition}\textbf{Piecewise harmonic function of order~$m\,\in\,\N^\star$}\\

\noindent  Given a strictly positive integer~$m$, a real valued function~$u$, defined on the set of vertices of~$\Gamma_{\cal W}$, is said to be \textbf{piecewise harmonic function of order~$m$} if, for any word~${\cal M}$ of length~$ m$,~$u\circ T_{\cal M}$ is harmonic of order~$m$.

\end{definition}

\vskip 1cm

\begin{definition}\textbf{Existence domain of the Laplacian, for a continuous function on the graph~$\Gamma_{\cal W}$ } (see \cite{Beurling1985})\\

\label{Lapl}
\noindent We will denote by~$\text{dom}\, \Delta$ the existence domain of the Laplacian, on the graph~$\Gamma_{\cal W}$, as the set of functions~$u$ of~$\text{dom}\, \mathcal{E}$ such that there exists a continuous function on~$\Gamma_{\cal W}$, denoted~$\Delta \,u$, that we will call \textbf{Laplacian of~$u$}, such that :
$$\mathcal{E}(u,v)=-\displaystyle \int_{{\cal D} \left ( \Gamma_{\cal W} \right)} v\, \Delta u   \,d\mu \quad \text{for any } v \,\in \,\text{dom}_0\, \mathcal{E}$$
\end{definition}

\vskip 1cm

\begin{definition}\textbf{Harmonic function}\\

\noindent A function~$u$ belonging to~\mbox{$\text{dom}\,\Delta$} will be said to be \textbf{harmonic} if its Laplacian is equal to zero.
\end{definition}

\vskip 1cm

\begin{notation}

In the following, we will denote by~${\cal H}_0\subset \text{dom}\, \Delta$ the space of harmonic functions, i.e. the space of functions~$u \,\in\,\ \text{dom}\, \Delta$ such that:

$$\Delta\,u=0$$

\noindent Given a natural integer~$m$, we will denote by~${\cal S} \left ({\cal H}_0,V_m \right)$ the space, of dimension~$N_b^m$, of spline functions " of level~$m$", ~$u$, defined on~$\Gamma_{\cal W}$, continuous, such that, for any word~$\cal M$ of length~$m$,~\mbox{$u \circ T_{\cal M}$} is harmonic, i.e.:

$$\Delta_m \, \left ( u \circ T_{\cal M} \right)=0$$

\end{notation}

\vskip 1cm

\begin{pte}

For any natural integer~$m$:

$${\cal S} \left ({\cal H}_0,V_m \right )\subset  \text{dom }{ \cal E}$$

\end{pte}
\vskip 1cm

\subsection{Explicit determination of the Laplacian of a function~$u$ of~\mbox{$\text{dom}\,\Delta$}}

\begin{definition}\textbf{Spline functions on~${\cal D}\left (\Gamma_{{\cal W}_m}\right)$,~$m\,\in\,\N^\star$}\\



\noindent Given a strictly positive integer~$m$, let us consider a vertex~$X$ of the graph~${\Gamma}_{{\cal W}_m}$.  Two configurations can occur:\\

\begin{enumerate}
	
	\item[\emph{i}.] the vertex~$X$ belongs to one and only one polygon with~$N_b$ sides,~${\cal P}_{m,j}$,~\mbox{$0 \leq j \leq N_b^m-1$}.\\
	
	In this case, if one considers the spline functions~$\psi_{Z}^{m}$ which correspond to the~$N_b-1$ vertices of this polygon distinct from~$X$:
	
	$$\displaystyle \sum_{Z \, \text{vertex of~${\cal P}_{m,j}$}}  \displaystyle\int_{{\cal D}\left (\Gamma_{\cal W}\right)} \psi_{Z}^{m}\,d\mu =\mu\left ({\cal P}_{m,j} \right)$$
	
	\noindent i.e., by symmetry:
	$$N_b \,  \displaystyle\int_{{\cal D}\left (\Gamma_{\cal W}\right)} \psi_{X}^{m}\,d\mu = \mu\left ({\cal P}_{m,j} \right)$$
	
	\noindent Thus:

	$$  \displaystyle\int_{{\cal D}\left (\Gamma_{\cal W}\right)}  \psi_{X}^{m}\,d\mu =\displaystyle\frac{1}{N_b}\, \mu\left ({\cal P}_{m,j} \right)= \frac{{\cal A}_m}{N_b}$$

\noindent where we retrieve the quantities~${\cal A}_m$ of definition~\ref{Am}.

	\begin{figure}[h!]
		\center{\psfig{height=8cm,width=10cm,angle=0,file=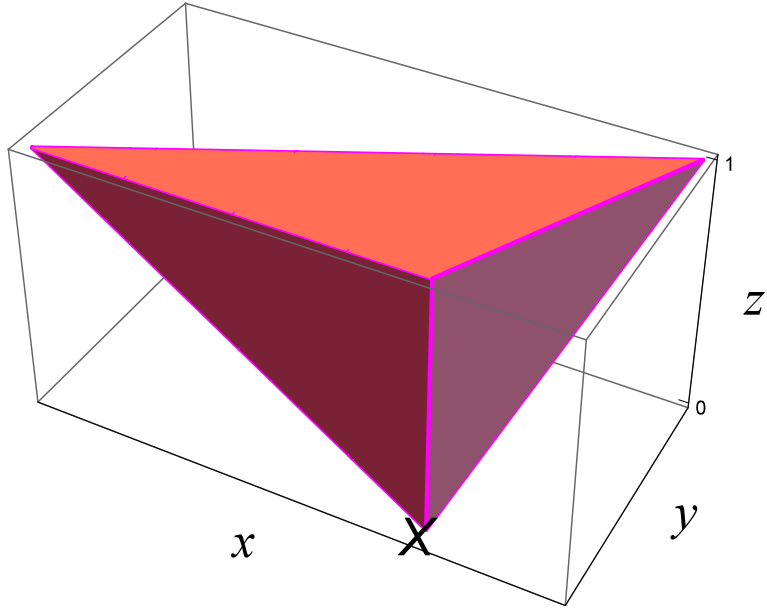}}\\
		\caption{The graph of a spline function~$\psi_{X}^{m}$,~$m\,\in\,\N$, in the case~$N_b=3$.}
	\end{figure}
	
	\vskip 1cm

	\item[\emph{ii}.] the vertex~$X$ is the intersection point of two polygons with~$N_b$ sides,~${\cal P}_{m,j}$ and~${\cal P}_{m,j+1}$,~\mbox{$0 \leq j \leq N_b^m-2$}.\\
	
	\noindent On has then to take into account the contributions of both polygons, which leads to:

	$$   \displaystyle\int_{{\cal D}\left (\Gamma_{\cal W}\right)}  \psi_{X}^{m}\,d\mu =\displaystyle\frac{1}{2\,N_b}\, \left \lbrace \mu\left ({\cal P}_{m,j} \right)+\mu\left ({\cal P}_{m,j+1} \right) \right\rbrace = \frac{{\cal A}_m}{N_b}$$

	\noindent where, again, we retrieve the quantities~${\cal A}_m$ of definition~\ref{Am}.

\end{enumerate}

\end{definition}
\vskip 1cm

\begin{remark} 

\noindent As it is explained in~\cite{StrichartzLivre2006}, one has just to reason by analogy with the dimension~1, more particulary, the unit interval~$I=[0,1]$, of extremities~$X_0=(0,0)$, and~$X_1=(1,0)$. The functions~$\psi_{X_1}$ and~$\psi_{X_2}$ such that, for any~$Y$ of~$\R^2$ :

$$\psi_{X_1} (Y)=\delta_{X_1Y} \quad  ,  \quad \psi_{X_2} (Y)=\delta_{X_2Y}   $$

\noindent are, in the most simple way, tent functions. For the standard measure, one gets values that do not depend on~$X_1$, or~$X_2$ (one could, also, choose to fix~$X_1$ and~$X_2$ in the interior of~$I$) :

$$\displaystyle\int_{I}  \psi_{X_1}\, d\mu =\displaystyle\int_{I}  \psi_{X_2}\, d\mu=\displaystyle \frac{1}{2}$$

\noindent (which corresponds to the surfaces of the two tent triangles.) \\

\begin{figure}[h!]
	\center{\psfig{height=8cm,width=10cm,angle=0,file=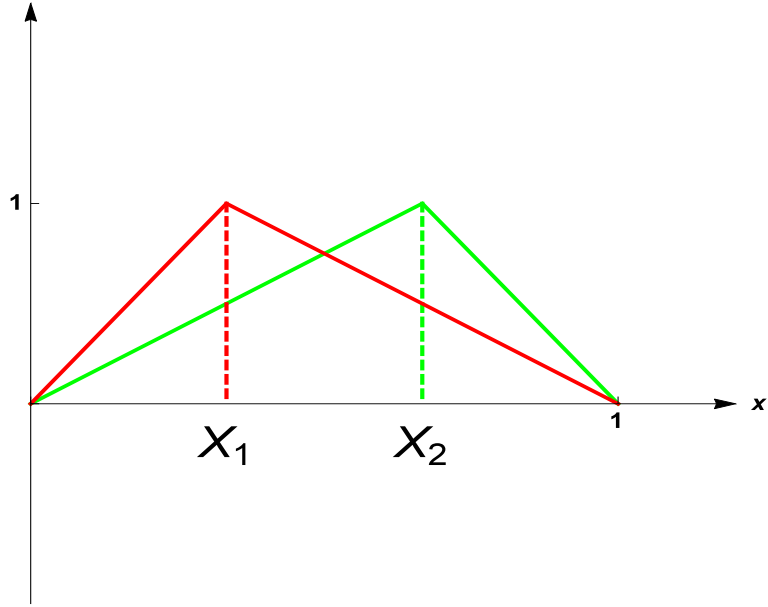}}\\
	\caption{The graphs of the spline functions~$\psi_{X_1}$ and~$\psi_{X_2}$.}
\end{figure}

\noindent In our case, we have just build the pendant, by no longer reasoning on the unit interval, but on our $N_b$-gons. \\
 
\end{remark}

\vskip 1cm

\begin{pte}
Let~$m$ be a strictly positive integer,~$X \,\notin\,V_0$ a vertex of the graph~$\Gamma_{\cal W}$, and~\mbox{$\psi_X^{m}\,\in\,{\cal S} \left ({\cal H}_0,V_m \right)$} a spline  function such that:

$$\psi_X^{m}(Y)=\left \lbrace \begin{array}{ccc}\delta_{XY} & \forall& Y\,\in \,V_m \\
 0 & \forall& Y\,\notin \,V_m \end{array} \right. \quad,  \quad \text{where} \quad    \delta_{XY} =\left \lbrace \begin{array}{ccc}1& \text{if} & X=Y\\ 0& \text{else} &  \end{array} \right.$$

\noindent Then, since~$X\, \notin \,V_0$: $\psi_X^{m}\,\in \,\text{dom}_0\, \mathcal{E}$.

\noindent Let us first note that:

\footnotesize

$$\begin{array}{ccc} {\cal E}  \left ( u ,\psi_X^{m} \right) &=&\displaystyle\sum_{Y \in V_m,\,Y\underset{m}{\sim} X} \Delta_m u(X) \, \left ( \psi_X^{m}(Y)(Y)- \psi_X^{m}(Y) \right) \, \left (  \displaystyle\int_{{\cal D} \left ( \Gamma_{\cal W} \right)}  \psi_X^{m}\, d\mu  \right )\\
&= &  \displaystyle \frac{ N_b^m}{h_m^2}\, \displaystyle\sum_{Y \in V_m,\,Y\underset{m}{\sim} X}
\left (u(Y)-u(X)\right)\,  \left (  \displaystyle\int_{{\cal D} \left ( \Gamma_{\cal W} \right)}  \psi_X^{m}\, d\mu  \right) \\
&= &\displaystyle \frac{ N_b^m}{h_m^2}\, \displaystyle\int_{{\cal D} \left ( \Gamma_{\cal W} \right)} 
\left (u(Y)-u(X)\right)\,  \psi_X^{m}(Y)\, d\mu \  \\
\end{array}$$

\normalsize

\noindent For any function~$u$ of~$\text{dom}\, \mathcal{E}$, such that its Laplacian exists, definition (\ref{Lapl}) applied to~$\psi_X^{m}$ leads to: 

$$\mathcal{E}(u,\psi_X^{m})= -N_b^{ m}\,\Delta_m u(X)\, \left (  \displaystyle\int_{{\cal D} \left ( \Gamma_{\cal W} \right)}  \psi_X^{m}\, d\mu  \right )=- \displaystyle\int_{{\cal D}\left (\Gamma_{\cal W}\right)}  \psi_X^{m}\,\Delta u  \, d\mu \approx -\Delta  u(X)\, \displaystyle\int_{{\cal D} \left ( \Gamma_{\cal W} \right)}  \psi_X^{m}\, d\mu$$

\noindent since~$\Delta u$ is continuous on~$ \Gamma_{\cal W} $, and the support of the spline function~$\psi_X^{m}$ is close to~$X$:

$$\displaystyle\int_{{\cal D} \left ( \Gamma_{\cal W} \right)}  \psi_X^{m}\,\Delta u  \, d\mu \approx -\Delta  u(X)\, \displaystyle\int_{{\cal D} \left ( \Gamma_{\cal W} \right)}  \psi_X^{m}\, d\mu$$

\noindent By passing through the limit when the integer~$m$ tends towards infinity, one gets:

$$ \displaystyle \lim_{m \to + \infty} \displaystyle\int_{{\cal D} \left ( \Gamma_{\cal W} \right)}  \psi_X^{m}\,\Delta_m u  \, d\mu=
 \Delta  u(X)\,\displaystyle \lim_{m \to + \infty}   \displaystyle\int_{{\cal D} \left ( \Gamma_{\cal W} \right)}  \psi_X^{m}\, d\mu$$

\noindent i.e.:

$$  \Delta  u(X)= \displaystyle \lim_{m \to + \infty}  N_b^m\,\Delta_m u(X)\,$$


\end{pte}

\vskip 1cm

\vskip 1cm
\begin{theorem}

\noindent Let~$u$ be in~\mbox{$\text{dom}\,\Delta$}. Then, the sequence of functions~$\left (f_m \right)_{m\in\N}$ such that, for any natural integer~$m$, and any~$X$ of~\mbox{$V_\star\setminus V_0$} :

 $$f_m(X)=N_b^{ m}\,\Delta_m \,u(X) $$

 \noindent  converges uniformly towards~$\Delta\,u$, and, reciprocally, if the sequence of functions~$\left (f_m \right)_{m\in\N}$ converges uniformly towards a continuous function on~\mbox{$V_\star\setminus V_0$}, then:

 $$u \,\in\, \text{dom}\,\Delta$$
\end{theorem}

\vskip 1cm
\begin{proof}

\noindent Let~$u$ be in~\mbox{$\text{dom}\,\Delta$}.  Since~$u$ belongs to~\mbox{$\text{dom}\,\Delta$}, its Laplacian~$\Delta\,u$ exists, and is continuous on the graph~$\Gamma_{\cal W}$. The uniform convergence of the sequence~$\left (f_m \right)_{m\in\N}$ follows.\\

\noindent Reciprocally, if the sequence of functions~$\left (f_m \right)_{m\in\N}$ converges uniformly towards a continuous function on~\mbox{$V_\star\setminus V_0$}, the, for any natural integer~$m$, and any~$v$ belonging to~\mbox{$\text{dom}_0\,{\cal E}$}:

$$\begin{array}{ccc} {\cal{E}}_{m }(u,v)
  &=&  \displaystyle \sum_{(X,Y) \,\in \, V_m^2,\,X  \underset{m }{\sim}  Y} 
  \frac{N_b^m}{h_m^2}\,\left (u_{\mid V_m}(X)-u_{\mid V_m}(Y)\right )\,\left(v_{\mid V_m}(X)-v_{\mid V_m}(Y)\right) \\
  &=& \displaystyle \sum_{(X,Y) \,\in \, V_m^2,\,X  \underset{m }{\sim}  Y} \frac{N_b^m}{h_m^2}\,\left (u_{\mid V_m}(Y)-u_{\mid V_m}(X )\right )\,\left(v_{\mid V_m}(Y)-v_{\mid V_m}(X)\right) \\
  &=&- \displaystyle \sum_{X \,\in \,V_m\setminus V_0 } \frac{N_b^m}{h_m^2}\,\sum_{Y\,\in \,V_m, \, Y  \underset{m }{\sim}  X} v_{\mid V_m}(X)\,\left (u_{\mid V_m}(Y)-u_{\mid V_m}(X)\right )   \\
  & &- \displaystyle \sum_{X \,\in \,  V_0 } \frac{N_b^m}{h_m^2}\,\sum_{Y\,\in \,V_m ,\, Y  \underset{m }{\sim}  X} v_{\mid V_m}(X)\,\left (u_{\mid V_m}(Y)-u_{\mid V_m}(X)\right )   \\
  &=&- \displaystyle \sum_{X \,\in \,V_m\setminus V_0 } N_b^m\,v(X)\,\Delta_m \,u(X)  \\
    &=&- \displaystyle \sum_{X \,\in \,V_m\setminus V_0 } v(X)    \, N_b^m\,  \Delta_m \,u(X)  \\
  \end{array}
  $$

\noindent Let us note that any~$X$ of~$V_m\setminus V_0$ admits exactly two adjacent vertices which belong to~$V_m\setminus V_0$, which accounts for the fact that the sum

 $$\displaystyle \sum_{X \,\in \,V_m\setminus V_0 } \frac{N_b^m}{h_m^2}\,\sum_{Y\,\in \,V_m\setminus V_0 ,\, Y  \underset{m }{\sim}  X} v(X)\,\left (u_{\mid V_m}(Y)-u_{\mid V_m}(X)\right)$$
 \noindent has the same number of terms as:

 $$ \displaystyle \sum_{(X,Y) \,\in \,(V_m\setminus V_0)^2,\,X  \underset{m }{\sim}  Y} \frac{N_b^m}{h_m^2}\,\left (u_{\mid V_m}(Y)-u_{\mid V_m}(X)\right )\,\left(v_{\mid V_m}(Y)-v_{\mid V_m}(X)\right)  $$

 \noindent For any natural integer~$m$, we introduce the sequence of functions~$\left (f_m \right)_{m\in\N}$ such that, for any~$X$ of~$V_m\setminus V_0$:

 $$f_m(X)=N_b^m\,\Delta_m \,u(X) $$

 \noindent The sequence~$\left (f_m \right)_{m\in\N}$ converges uniformly towards~$\Delta\,u$. 


\end{proof}


\section{Normal derivatives}

Let us go back to the case of a function~$u$ twice differentiable on~$I=[0,1]$, that does not vanish in~0 and~1:
$$  \displaystyle \int_0^1 \left (\Delta u \right)(x)\,v(x)\,dx=-   \displaystyle \int_0^1 u'(x)\,v'(x)\,dx+ u '(1)\,v (1)-u '(0)\,v (0) $$

\noindent The normal derivatives:

$$\partial_n u(1)=u'(1) \quad \text{and} \quad \partial_n u(0)=u'(0) $$

\noindent appear in a natural way. This leads to:

$$  \displaystyle \int_0^1 \left (\Delta u \right)(x)\,v(x)\,dx=-   \displaystyle \int_0^1 u'(x)\,v'(x)\,dx+ \displaystyle \sum_{\partial\, [0,1]} v\,\partial_n\,u $$

\noindent One meets thus a particular case of the Gauss-Green formula, for an open set~$\Omega$ of~$\R^d$,~$d \,\in\,\N^\star$:

$$  \displaystyle \int_\Omega \nabla\,  u\, \nabla \, v \, d \mu= -\displaystyle \int_\Omega  \left (\Delta u \right) \,v \,d\mu + \displaystyle \int_{\partial\, \Omega } v\,\partial_n\,u \,d\sigma$$

\noindent where~$\mu$ is a measure on~$\Omega $, and where~$d\sigma$ denotes the elementary surface on~$\partial\, \Omega $.\\

\noindent In order to obtain an equivalent formulation in the case of the graph~$\Gamma_{\cal W}$, one should have, for a pair of functions~$(u,v)$ continuous on~$\Gamma_{\cal W}$ such that~$u$ has a normal derivative:

$$  {\cal E}(u,v)= -\displaystyle \int_\Omega  \left (\Delta u \right) \,v \,d\mu + \displaystyle \sum_{V_0} v\,\partial_n\,u  $$

\noindent For any natural integer~$m$ :

$$\begin{array}{ccc}
 {\cal{E}}_{m }(u,v)
  &=&  \displaystyle \sum_{(X,Y) \,\in \, V_m^2,\,X  \underset{m }{\sim}  Y} \frac{N_b^m}{h_m^2}\,\left (u_{\mid V_m}(Y)-u_{\mid V_m}(X)\right )\,\left(v_{\mid V_m}(Y)-v_{\mid V_m}(X)\right) \\
  &=&- \displaystyle \sum_{X \,\in \,V_m\setminus V_0 } \frac{N_b^m}{h_m^2}\,\sum_{Y\,\in \,V_m  ,\, Y  \underset{m }{\sim}  X} v_{\mid V_m}(X)\,\left (u_{\mid V_m}(Y)-u_{\mid V_m}(X)\right )   \\
   &&-  \displaystyle \sum_{X \,\in \,  V_0 } \frac{N_b^m}{h_m^2}\,\sum_{Y\,\in \,  V_m ,\, Y  \underset{m }{\sim}  X} v_{\mid V_m}(X)\,\left (u_{\mid V_m}(Y)-u_{\mid V_m}(X)\right )  \\
  &=&- \displaystyle \sum_{X \,\in \,V_m\setminus V_0 } v_{\mid V_m}(X)\,N_b^m\,\Delta_m \,u_{\mid V_m}(X)  \\
   &&+ \displaystyle \sum_{X \,\in \,  V_0 } \sum_{Y\,\in \,  V_m ,\, Y  \underset{m }{\sim}  X} \frac{N_b^m}{h_m^2}\,v_{\mid V_m}(X)\,\left (u_{\mid V_m}(X)-u_{\mid V_m}(Y)\right )  \\
  \end{array}
  $$




 \noindent We thus come across an analogous formula of the Gauss-Green one, where the role of the normal derivative is played by:

 $$ \displaystyle \sum_{X \,\in \,  V_0 } \frac{N_b^m}{h_m^2}\,\sum_{Y\,\in \,  V_m ,\, Y  \underset{m }{\sim}  X} \,\left (u_{\mid V_m}(X)-u_{\mid V_m}(Y)\right ) \,\frac{{\cal A}_m}{N_b}$$

 \vskip 1cm

 \begin{definition}

 \noindent For any~$X$ of~$V_0$, and any continuous function~$u$ on~$\Gamma_{\cal W}$, we will say that~$u$ admits a normal derivative in~$X$, denoted by~$\partial_n\,u(X)$, if:

 $$ \displaystyle \lim_{m \to + \infty}  \frac{N_b^m}{h_m^2}\,\sum_{Y\,\in \,  V_m ,\, Y  \underset{m }{\sim}  X} \,\left (u_{\mid V_m}(X)-u_{\mid V_m}(Y)\right ) \,\frac{{\cal A}_m}{N_b}< + \infty $$

 \noindent We will set:

 $$\partial_n\,u(X) = \displaystyle \lim_{m \to + \infty}  \frac{N_b^m}{h_m^2}\,\sum_{Y\,\in \,  V_m ,\, Y  \underset{m }{\sim}  X} \,\left (u_{\mid V_m}(X)-u_{\mid V_m}(Y)\right ) \,\frac{{\cal A}_m}{N_b}< + \infty $$

 \end{definition}

 \vskip 1cm

 \begin{definition}

\noindent For any natural integer~$m$, any~$X$ of~$V_m$, and any continuous function~$u$ on~$\Gamma_{\cal W}$, we will say that~$u$ admits a normal derivative in~$X$, denoted by~$\partial_n\,u(X)$, if:
 $$ \displaystyle \lim_{k \to + \infty}  \frac{r^{-k}}{h_k^2}\,\sum_{Y\,\in \,  V_k ,\, Y  \underset{k }{\sim}  X} \,\left (u_{\mid V_k} (X)-u_{\mid V_k}(Y)\right )\,\frac{{\cal A}_k}{N_b} < + \infty $$

 \noindent We will set:

 $$\partial_n\,u(X) = \displaystyle \lim_{k \to + \infty}\frac{r^{-k}}{h_k^2}\,\sum_{Y\,\in \,  V_k ,\, Y  \underset{k }{\sim}  X} \,\left (u_{\mid V_k}(X)-u_{\mid V_k}(Y)\right ) \,\frac{{\cal A}_k}{N_b}< + \infty $$

 \end{definition}

 \vskip 1cm
 \begin{remark} One can thus extend the definition of the normal derivative of~$u$ to~$\Gamma_{\cal W}$.
 \end{remark}

 \vskip 1cm

 \begin{theorem}

\noindent Let~$u$ be in~\mbox{$\text{dom}\,\Delta$}. The, for any~$X$ of~$\Gamma_{\cal W}$,~$\partial_n\,u(X)$ exists. Moreover, for any~$v$ of~\mbox{$\text{dom}\,\cal E$}, et any natural integer~$m$, the Gauss-Green formula writes:

$$  {\cal E}(u,v) = -\displaystyle \int_{{\cal D} \left ( \Gamma_{\cal W}\right) } \left (\Delta u \right) \,v \,d\mu + \displaystyle \sum_{V_0} v\,\partial_n\,u  $$
 \end{theorem}
 \vskip 1cm

\section{Spectrum of the Laplacian}

In the following, let~$u$ be in~$\text{dom}\, \Delta$. We will apply the \emph{\textbf{spectral decimation method}} developed by~R.~S.~Strichartz \cite{StrichartzLivre2006}, in the spirit of the works of M.~Fukushima et T.~Shima \cite{Fukushima1994}. In order to determine the eigenvalues of the Laplacian~$\Delta\, u$ built in the above, we concentrate first on the eigenvalues~$\left (-{\Lambda_m}\right)_{m\in\N}$ of the sequence of graph Laplacians~$\left (\Delta_m \,u\right)_{m\in\N}$, built on the discrete sequence of graphs~$\left (\Gamma_{{ \cal W}_m}\right)_{m\in\N}$. For any natural integer~$m$, the restrictions of the eigenfunctions of the continuous Laplacian~$\Delta\,u$ to the graph~$\Gamma_{{ \cal W}_m}$ are, also, eigenfunctions of the Laplacian~$\Delta_m$, which leads to recurrence relations between the eigenvalues of order~$m$ and~$m+1$.

\vskip 1cm

We thus aim at determining the solutions of the eigenvalue equation:

$$-\Delta\,u=\Lambda\,u \quad  \text{on }\Gamma_{{ \cal W} }$$

\noindent as normalized limits, when the integer~$m$ tends towards infinity, of the solutions of:

$$-\Delta_m\,u=\Lambda_m\,u \quad  \text{on }V_m\setminus V_0$$

\noindent Let~$m \geq 1$. We consider an eigenfunction~$u_{m-1}$ on~\mbox{$V_{m-1}\setminus V_0$}, for the eigenvalue~$\Lambda_{m-1}$. The aim is to extend~$u_{m-1}$ on~\mbox{$V_m\setminus V_0$} in a function~$u_m$, which will itself be an eigenfunction of~$\Delta_m$, for the eigenvalue~$\Lambda_m$, and, thus, to obtain a recurrence relation between the eigenvalues~$\Lambda_m$ and~$\Lambda_{m-1}$. Given three consecutive vertices of~$\Gamma_{{ \cal W}_{m-1}}$,~$X_k$,~$X_{k+1}$,~$X_{k+2}$, where~$k$ denotes a generic natural integer, we will denote by~$Y_{k+1}$,~$\hdots$,~$Y_{k+N_b-1}$ the points of~\mbox{$V_m\setminus V_{m-1}$} such that:~$Y_{k+1}$,~$\hdots$,~$Y_{k+N_b-1}$ are  between~$X_k$ and~$X_{k+1}$, and by~$Y_{k+N_b+1}$,~$\hdots$,~$Y_{k+2\,N_b-1}$, the points of~\mbox{$V_m\setminus V_{m-1}$} such that:~$Z_{k+1}$,~$\hdots$,~$Z_{k+N_b-1}$ are between~$X_{k+1}$ and~$X_{k+2}$. For the sake of consistency, let us set:

$$Y_{k+N_b }=X_{k+1} \quad \text{and} \quad  Y_{k+2\,N_b }=X_{k+2}$$

 \begin{figure}[h!]
 \center{\psfig{height=8cm,width=14cm,angle=0,file=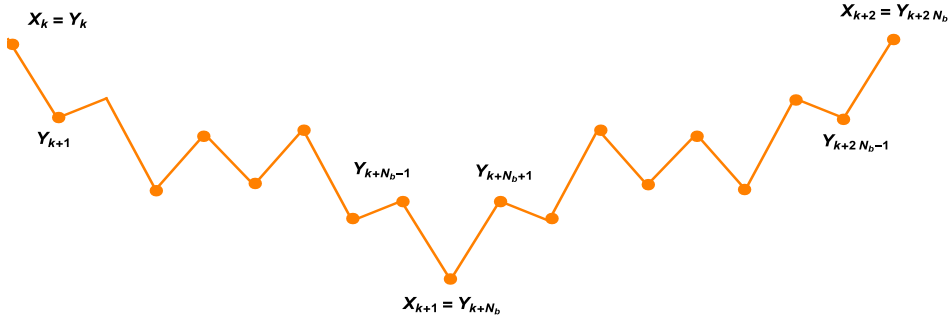}}\\
\caption{The points~$X_k$,~$X_{k+1}$,~$X_{k+2}$, and~$Y_{k }$,~$\hdots$,~$Y_{k+N_b }$,~$\hdots$,~$Y_{k+2\,N_b }$.}
 \end{figure}

\vskip 1cm
\newpage



\vskip 1cm

\noindent The values of~$u_{m-1}$ in~$X_k$,~$X_{k+1}$,~$X_{k+2}$ are thus supposed to be known. \\
\noindent The eigenvalue equation in~$\Lambda_m$ leads to the following systems:

$$\left \lbrace
\begin{array}{ccccc}
   \left \lbrace \Lambda_{m }-2\, h_{m-1}^2 \right\rbrace \,u_m\left ( Y_{k+1} \right) &=& -h_{m-1}^2\, u_{m-1}\left ( X_k \right) -h_{m-1}^2\,u_m\left ( Y_{k+2} \right) &&\\
   \left \lbrace \Lambda_{m }-2\, h_{m-1}^2\, \right\rbrace \,u_m\left ( Y_{k+i} \right) &=& - h_{m-1}^2\,u_{m }\left ( Y_{k+i-1} \right) -h_{m-1}^2\,u_m\left ( Y_{k+i+1} \right)
& , & 1 \leq i \leq N_b-3 \\
  \left \lbrace \Lambda_{m }-2\, h_{m-1}^2\, \right\rbrace \,u_m\left ( Y_{k+N_b-1} \right) &=& - h_{m-1}^2\,u_{m -1}\left ( X_{k+1} \right) -h_{m-1}^2\,u_{m }\left ( Y_{k +N_b- 2} \right) &&
\end{array}
\right.$$

\noindent and :

$$\left \lbrace
\begin{array}{ccccc}   \left \lbrace \Lambda_{m }-2 \right\rbrace \,u_m\left ( Y_{k+N_b+1} \right) &=& - h_{m-1}^2\,u_{m-1}\left ( X_{k+1} \right) -h_{m-1}^2\,u_m\left ( Y_{k+N_b+2}\right) &&\\
 \left \lbrace \Lambda_{m }-2\,h_{m-1}^2\, \right\rbrace \,u_m\left ( Y_{k+N_b+i} \right) &=& - h_{m-1}^2\,u_{m }\left (Y_{k+N_b+i-1} \right) -h_{m-1}^2\,u_m\left ( Y_{k+N_b+i+1} \right) &,&  1 \leq i \leq N_b-3 \\
   \left \lbrace \Lambda_{m }-2\,  h_{m-1}^2\,\right\rbrace \,u_m\left ( Y_{k+2\,N_b-1} \right) &=& - h_{m-1}^2\,u_{m -1}\left ( X_{k+2} \right) -h_{m-1}^2\,u_{m }\left ( Y_{k+2\,N_b-2} \right) &&\end{array}
\right.$$

\noindent For the sake of simplicity, we set:

$$ \Lambda_{m }= h_{m-1}^2\,  \tilde{\Lambda}_{m }$$

\noindent The sequence~$\left ( u_m\left ( Y_{k+i} \right) \right)_{0 \leq i \leq 2\,N_b}$ satisfies a second order recurrence relation, the characteristic equation of which is:

 $$ r^2+  \left \lbrace\tilde{\Lambda}_{m }-2 \right\rbrace \,r+1=0 $$

 \noindent The discriminant is:

$$ \delta_m=   \left \lbrace \tilde{\Lambda}_{m }-2 \right\rbrace ^2-4= \omega_m^2 \quad, \quad \omega_m\,\in\,\C$$





\noindent The roots~$r_{1,m}$ and~$r_{2,m}$ of the characteristic equation are the scalar given by:

$$r_{1,m}=\displaystyle \frac{ 2-\tilde{\Lambda}_{m }   -\omega_m}{2}  \quad, \quad r_{2,m}=\displaystyle \frac{  2-\tilde{\Lambda}_{m } +\omega_m}{2}$$

\noindent One has then, for any natural integer~$i$ of~\mbox{$\left \lbrace 0,\hdots,2\,N_b \right \rbrace $} :

$$ u_m\left ( Y_{k+i } \right) = \alpha_m\, r_{1,m}^i +\beta_m\, r_{2,m}^i$$

\noindent where~$\alpha_m$ and~$\beta_m$ denote scalar constants.

\noindent The extension~$u_m$ of~$u_{m-1}$ to~\mbox{$V_m\setminus V_0$} has to be an eigenfunction of~$\Delta_m$, for the eigenvalue~$\Lambda_m$.\\
\noindent Since~$u_{m-1}$ is an eigenfunction of~$\Delta_{m-1}$, for the eigenvalue~$\Lambda_{m-1}$, the sequence~$\left ( u_{m-1}\left ( X_{k+i} \right) \right)_{0 \leq i \leq  N_b}$ must itself satisfy a second order linear recurrence relation which be the pendant, at order~$m$, of the one satisfied by the sequence~$\left ( u_m\left ( Y_{k+i} \right) \right)_{0 \leq i \leq 2\,N_b}$, the characteristic equation of which is:

 $$ h_{m-2}^2\, \left \lbrace \Lambda_{m-1 }-2 \right\rbrace \,r= - 1 -r^2$$

\noindent and discriminant:

$$ \delta_{m-1}=  \left \lbrace \Lambda_{m-1 }-2 \right\rbrace ^2-4= \omega_{m-1}^2 \quad, \quad \omega_{m-1}\,\in\,\C$$

\noindent The roots~$r_{1,m-1}$ and~$r_{2,m-1}$ of this characteristic equation are the scalar given by:

$$r_{1,m-1}=\displaystyle \frac{ 2-\tilde{\Lambda}_{m-1} -\omega_{m-1}}{2} \quad, \quad r_{2,m-1}=\displaystyle \frac{ 2-\tilde{\Lambda}_{m-1} +\omega_{m-1}}{2}$$

\noindent For any integer~$i$ of~\mbox{$\left \lbrace 0,\hdots, N_b \right \rbrace $}:

$$ u_{m-1}\left ( Y_{k+i } \right) = \alpha_{m-1}\, r_{1,m-1}^i +\beta_{m-1}\, r_{2,m-1}^i$$

\noindent where~$\alpha_{m-1}$ and~$\beta_{m-1}$ denote scalar constants.

\noindent From this point, the compatibility conditions, imposed by spectral decimation, have to be satisfied:

$$  \left \lbrace \begin{array}{ccc}
u_{m }\left ( Y_{k  } \right)&=&u_{m-1}\left ( X_{k  } \right)  \\
u_{m }\left ( Y_{k+N_b  } \right)&=&u_{m-1}\left ( X_{k+1 } \right)  \\
u_{m }\left ( Y_{k+2\,N_b } \right)&=&u_{m-1}\left ( X_{k+2 } \right)  \\
\end{array}\right.$$

\noindent i.e.:

$$  \left \lbrace \begin{array}{ccccc}
 \alpha_{m }  +\beta_{m } &=& \alpha_{m-1}  +\beta_{m-1}  & {\cal C}_{ m}\\
 \alpha_{m }\, r_{1,m }^{N_b} +\beta_{m }\, r_{2,m }^{N_b}&=& \alpha_{m-1}\, r_{1,m-1}  +\beta_{m-1}\, r_{2,m-1} & {\cal C}_{1, m} \\
 \alpha_{m }\, r_{1,m }^{2\,N_b} +\beta_{m }\, r_{2,m }^{2\,N_b}&=& \alpha_{m-1}\, r_{1,m-1}^2 +\beta_{m-1}\, r_{2,m-1}^2 & {\cal C}_{2, m}\\
\end{array}\right. $$

  \noindent where, for any natural integer~$m$,~$\alpha_{m }$ and~$\beta_{m }$ are scalar constants (real or complex).\\

\noindent Since the graph~$\Gamma_{{\cal W}_{ m-1}}$ is linked to the graph~$\Gamma_{{\cal W}_{ m }}$ by a similar process to the one that links~$\Gamma_{{\cal W}_{  1}}$ to~$\Gamma_{{\cal W}_{ 0}}$, one can legitimately consider that the constants~$\alpha_m$ and~$\beta_m$ do not depend on the integer~$m$:

$$\forall\,m\,\in\,\N^\star \, : \quad \alpha_m=\alpha \,\in\,\R  \quad, \quad \beta_m=\beta \,\in\,\R $$

\noindent The above system writes:

$$  \left \lbrace \begin{array}{ccc}
 \alpha \, r_{1,m }^{N_b} +\beta \, r_{2,m }^{N_b}&=& \alpha \, r_{1,m-1}  +\beta \, r_{2,m-1}  \\
 \alpha \, r_{1,m }^{2\,N_b} +\beta \, r_{2,m }^{2\,N_b}&=& \alpha \, r_{1,m-1}^2 +\beta \, r_{2,m-1}^2 \\
\end{array}\right.$$

\noindent One has then to consider the following configurations:
\begin{enumerate}
\item[\emph{i}.] \underline{First case:}  \\

For any natural integer~$m$ :

$$r_{1,m }\,\in\,\R \quad , \quad r_{2,m }\,\in\,\R$$

\noindent and, more precisely:

$$r_{1,m }<0 \quad , \quad r_{2,m }<0$$

\noindent since the function~$\varphi$, which, to any real number~$x \geq 4$, associates:

$$\varphi(x)=  \displaystyle \frac{2-x+\varepsilon \, \sqrt{\left \lbrace x-2 \right\rbrace ^2-4}}{2}
\quad , \quad \varepsilon \,\in \, \left \lbrace -1, 1 \right \rbrace $$



\noindent is strictly increasing on~$]4,+ \infty[$. Due to its continuity, is is a bijection of~$]4,+ \infty[$ on~\mbox{$ \varphi\left (  ]4,+ \infty[  \right)=]  -1,0 [ $}.\\

This configuration only occurs in the cases when the natural integer~$N_b$ is an odd number.
\noindent Let us introduce the function~$\phi$, which, to any real number~~$x \geq 2$, associates:

$$\phi(x)= | \varphi(x)|= \displaystyle \frac{-2+x-\varepsilon \,\sqrt{\left \lbrace x-2 \right\rbrace ^2-4}}{2}$$

\noindent where~\mbox{$\varepsilon\,\in\,\left \lbrace -1, 1 \right \rbrace$}.

\noindent The function~$\phi$ is a bijection of~$]4,+ \infty[$ on~\mbox{$ \phi\left (  ]4,+ \infty[  \right)=] 0,1 [ $}.
We will denote by~$\phi^{-1}$ its inverse bijection:

$$\forall\, \,x \,\in\,]0,1[ \, : \quad \phi^{-1}(x)= \displaystyle \frac{(y+1)^2}{y}$$. \\

\noindent One has then:

$$\varphi\left ( \tilde{\Lambda}_{m-1} \right) =\displaystyle \frac{2-\tilde{\Lambda}_{m-1}+\varepsilon\,\omega_{m-1}}{2}  \leq 0  $$

\noindent This yields:

$$
(-1)^{N_b} \, \left (\varphi\left ( \tilde{\Lambda}_{m } \right) \right)^{N_b} =   \varphi\left ( \tilde{\Lambda}_{m-1} \right) \leq 0  $$

\noindent which leads to:

$$   \phi\left (  \tilde{\Lambda}_{m } \right) = \left ( \phi\left ( \tilde{\Lambda}_{m-1} \right) \right)^{\frac{1}{N_b}}  $$






\noindent and:

$$\tilde{\Lambda}_{m } = \phi^{-1} \left (\left (\phi\left (\tilde{\Lambda}_{m-1 } \right ) \right )^{\frac{1}{N_b}} \right)
=\displaystyle \frac{\left \lbrace \left (\phi\left (\tilde{\Lambda}_{m-1 } \right ) \right )^{\frac{1}{N_b}}  +1 \right \rbrace^2}{\left (\phi\left (\tilde{\Lambda}_{m-1 } \right )
\right )^{\frac{1}{N_b}}   }
= \displaystyle \frac{\left \lbrace \left (\displaystyle \frac{-2+\tilde{\Lambda}_{m-1 }-\varepsilon \,\sqrt{\left \lbrace \tilde{\Lambda}_{m-1 }-2 \right\rbrace ^2-4}}{2}\right)
^{\frac{1}{N_b}}   +1 \right \rbrace^2}{\left (\displaystyle \frac{-2+\tilde{\Lambda}_{m-1 }-\varepsilon \,\sqrt{\left \lbrace \tilde{\Lambda}_{m-1 }-2 \right\rbrace ^2-4}}{2}\right)
^{\frac{1}{N_b}}   }
$$

\item[\emph{ii}.] \underline{Second case :}  \\

For any natural integer~$m$:

$$r_{1,m }\,\in\,\C\setminus \R \quad r_{2,m }= \overline{r_{1,m }}\,\in\,\C \setminus \R$$

\noindent Let us introduce:

$$\rho_{ m }= \left | r_{1,m }\right| \,\in\,\R^+ \quad , \quad \theta_m=\text{arg}\, r_{1,m } \quad \text{if} \quad r_{1,m } \neq 0$$

\noindent The above system writes:

$$  \left \lbrace \begin{array}{ccc}
\rho_{ m }^{ N_b}\, \left \lbrace \gamma \,\cos \left ( N_b\,\theta_m\right) +\delta\,  \sin \left ( N_b\,\theta_m\right) \right \rbrace &=& \rho_{ m-1 } \, \left \lbrace \gamma \,\cos \left (  \theta_{m-1}\right) +\delta\,  \sin \left (  \theta_{m-1}\right) \right \rbrace \\
\rho_{ m }^{ 2\,N_b}\, \left \lbrace \gamma \,\cos \left ( 2\, N_b\,\theta_m\right) +\delta\,  \sin \left (2\,N_b\,\theta_m\right) \right \rbrace &=&\rho_{ m-1 }^2 \, \left \lbrace \gamma \,\cos \left (  2\,\theta_{m-1}\right) +\delta\,  \sin \left ( 2\,\theta_{m-1}\right) \right \rbrace
\end{array}\right.$$

  \noindent where~$\gamma$ and~$\delta$ denote real constants.\\

  \noindent The system is satisfied if:

$$  \left \lbrace \begin{array}{ccc}
\rho_{ m }^{ N_b}  &=& \rho_{ m-1 }  \\
 \theta_m  &=&\displaystyle \frac{\theta_{m-1}}{N_b}
\end{array}\right.$$

\noindent and thus:

$$   \phi\left (  \tilde{\Lambda}_{m } \right) = \left ( \phi\left ( \tilde{\Lambda}_{m-1} \right) \right)^{\frac{1}{N_b}}  $$






\noindent which leads to the same relation as in the previous case:

$$ \tilde{\Lambda}_{m } = \phi^{-1} \left (\left (\phi\left (\tilde{\Lambda}_{m-1 } \right ) \right )^{\frac{1}{N_b}} \right)
=\displaystyle \frac{\left \lbrace \left (\phi\left (\tilde{\Lambda}_{m-1 } \right ) \right )^{\frac{1}{N_b}}  +1 \right \rbrace^2}{\left (\phi\left (\tilde{\Lambda}_{m-1 } \right )
\right )^{\frac{1}{N_b}}   }
= \displaystyle \frac{\left \lbrace \left (\displaystyle \frac{-2+\tilde{\Lambda}_{m-1 }-\varepsilon \,\sqrt{\left \lbrace \tilde{\Lambda}_{m-1 }-2 \right\rbrace ^2-4}}{2}\right)
^{\frac{1}{N_b}}   +1 \right \rbrace^2}{\left (\displaystyle \frac{-2+\tilde{\Lambda}_{m-1 }-\varepsilon \,\sqrt{\left \lbrace \tilde{\Lambda}_{m-1 }-2 \right\rbrace ^2-4}}{2}\right)
^{\frac{1}{N_b}}   }
$$

\noindent where~\mbox{$\varepsilon\,\in\,\left \lbrace -1, 1 \right \rbrace$}.

\end{enumerate}

\vskip 2cm

\centerline{\textbf{Thanks}}

\vskip 1cm
The author would like to thank JPG for his patient re-reading, and his very pertinent suggestions and advices, and Robert Strichartz, who suggested the introduction of specific energies to fully take into account the very specific geometry of the problem. Special thanks to GL, thanks to whom I searched in the right direction. All this helped improving the original work.

\vskip 2cm

\bibliographystyle{alpha}
\bibliography{BibliographieClaire}

\end{document}